\documentclass{article}[12pt]
\usepackage[english]{babel}
\usepackage{amsmath,amssymb,amsthm}
 \newtheorem{thm}{Theorem}[section]

 \newtheorem{prop}[thm]{Proposition}
 \theoremstyle{definition}
 \newtheorem{definition}[thm]{Definition}

  \newtheorem{example}[thm]{Example}
 \numberwithin{equation}{section}

\newcommand{\cc}[1]      {\overline{{#1}}}              

\newcommand{\Z}{\mathbb{Z}}

\newcommand{\intg}{\mathbb{Z}}
\newcommand{\bbnu}{ [[\nu]] }
\newcommand{\R}{\mathbb R}%
\newcommand*{\C}{\mathbb C}%

\newcommand{\g}{\mathfrak{g}}
\newcommand{\h}{\mathfrak{h}}

\newcommand{\SP}[1]      {\left\langle{#1}\right\rangle} 
\newcommand{\Bounded}[1][{}] {\operatorname{\mathfrak{B}}_{\scriptscriptstyle{#1}}}
\newcommand{\Anti}       {\Lambda}
\newcommand{\cech}{\v Cech}
\makeatletter
\def\operatorname#1{\mathop{\operator@font #1}\nolimits}%
\makeatother
\DeclareMathOperator{\ad}{ad}

\DeclareMathOperator{\Ad}{Ad}
\DeclareMathOperator{\id}{Id}
\DeclareMathOperator{\Aut}{Aut}
\DeclareMathOperator{\Ker}{Ker}
\DeclareMathOperator{\tr}{tr}

\DeclareMathOperator{\Der}{Der}

\DeclareMathOperator{\Hom}{Hom}
\DeclareMathOperator{\pr}{\mathsf{pr}}

\newcommand{\Image}{\operatorname{Im}}
\DeclareMathOperator{\Id}{Id}
\newcommand{\Cinftycf}   {\mathcal{C}^\infty_{\mathrm{cf}}}
\newcommand{\supp}       {\operatorname{\mathrm{supp}}}  
\newcommand{\e}{{\bf{e}}}

\newcommand*{\cyclic}{\mathop{\kern0.9ex{{+}\kern-2.2ex\raise-.29ex%
      \hbox{\Large\hbox{$\circlearrowright$}}}}\limits}
      
 \newcommand{\PO}{P}

\renewcommand{\H}{\mathcal{H}}
\newcommand{\half}{{\textstyle{\frac12}}}
\let\id\Id
\newcommand{\red}{{red}}
\newcommand{\Cinfty}{C^\infty}

\newcommand{\CiM}{C^\infty(M)}

\renewcommand{\L}{{\mathcal{L}}}
\renewcommand{\Im}{\operatorname{Im}}

\newcommand{\ehat[1]}{{\hat E}({L^{#1}})}
\newcommand{\CL}{{\mathcal C}_L}
\newcommand{\AntiC}      {\Lambda_{\mathbb{C}}}
\newcommand{\Sym}[1][\!]{\ensuremath{\mathcal{S}^{#1}}}
\newcommand{\Tens}[1][\!]{\ensuremath{\mathcal{T}^{#1}}}

\DeclareMathOperator{\un}{\mathbf{1}}

\newcommand{\deform}[1]  {\underline{\bf{#1}}}

\newcommand{\koszul}     {\partial}
\newcommand{\ins}        {\operatorname{\mathrm{i}}}
\newcommand{\prol}       {\operatorname{\mathrm{prol}}}
\newcommand{\nice}       {\mathrm{nice}}
\newcommand{\D}          {\operatorname{\mathrm{d}}} 
\newcommand{\Lie}        {\operatorname{\mathcal{L}\!}}    
\newcommand{\lie}[1]     {\mathfrak{#1}}
\newcommand{\image}      {\operatorname{\mathrm{im}}}
\newcommand{\zbar}{{\overline{z}}}
\newcommand{\wbar}{{\overline{w}}}
\newcommand{\Disk}{{\mathbb{D}}}
\def\cyclic{\mathop{\kern0.9ex{{+}
\kern-1.8ex\raise-.01ex\hbox{
{$\circlearrowright$}}}}\limits}

\newcommand{\starred}    {\mathbin{\star_\red}}

\newcommand{\starredk}   {\mathbin{\star_\red^{(\kappa)}}}

\newcommand{\deformiotak}{\deform{\iota}^{\deform{*}}_\kappa}
\newcommand{\bulletk}    {\mathbin{\bullet_\kappa}}

\title{Group Actions in Deformation Quantisation  }
\author{Simone Gutt\\[2mm]
{\footnotesize{  Acad\'emie Royale de Belgique}}\\
{\footnotesize{D\'{e}partement de Math\'{e}matique, Universit\'{e} Libre de Bruxelles}}\\
{\footnotesize{Campus Plaine, CP 218, Boulevard du Triomphe}}\\
 {\footnotesize{BE -- 1050 Bruxelles, Belgium}}\\
{\scriptsize{ Email:  sgutt@ulb.ac.be}}
}

\date{}

\begin{document}

\maketitle

\begin{abstract}

This set of notes corresponds to a mini-course given in September 2018 in Bedlewo; it does not contain any new result; 
it  complements -with intersection- the introduction  to
formal  deformation quantization and group actions published in \cite{bib:SG}, corresponding to a course given
in Villa de Leyva in July 2015.\\
After an introduction to the concept of deformation quantization, we briefly recall existence, classification and representation results for formal star products.
We come then to results concerning the notion of formal star products with symmetries; one has a Lie group action (or a Lie algebra action) compatible with the Poisson structure, and one wants to consider star products such that the Lie group acts by automorphisms (or the Lie algebra acts by derivations). We recall in particular the link between left invariant star products on Lie groups and Drinfeld twists, and the notion of universal deformation formulas.
Classically, symmetries are particularly interesting when they are implemented by a moment map and we give indications to build a corresponding quantum moment map. 
Reduction is a construction in classical mechanics with symmetries which allows to reduce the dimension of the manifold; we describe one of the various quantum analogues which  have been considered in the framework of formal deformation quantization.  
We end up by some considerations about convergence of star products.
 \end{abstract}

\section{Introduction to the notion of deformation quantization}

 A quantization gives a way to pass from a classical description
 to a quantum description of a physical system.  Since quantum theory provides a  description  of nature which is more fundamental than
 classical theory, one can wonder at the relevance of quantization.  Points in favour of such an attempt are the following: \\
 - Giving a priori a quantum description  of a physical system is  difficult, whilst  the classical description is often easier to obtain, so the classical description can be useful as a starting point  to find a  quantum description. \\
 - Any given physical theory remains  valid within a range of measurements, so that
any modified theory should 
give the same results in the initial range.\\ 
This second point  is an important  motivation of  the seminal idea of Moshe Flato that any new physical theory can appear as a deformation of the older one. 
In particular, the description of a system by classical mechanics
is good to describe the macroscopic non relativistic world.
Deformation Quantization was introduced by Flato,  Lichnerowicz and  Sternheimer  in \cite{bib:Flato} and developed in \cite{bib:Bayen}  to present quantum mechanics as a deformation of classical mechanics. 
One of the main feature of this quantization method is that  the emphasis is put on the algebra of observables.  They ``suggest that quantisation be understood as a deformation of the
structure of the algebra of classical observables rather than a radical change
in the nature of the observables.''

The classical description of mechanics
 in its  Hamiltonian formulation on
the motion space  (in general the quotient of the evolution space by the motion),  has for framework
a symplectic manifold $(M,\omega)$, or more generally a Poisson manifold $(M,\PO)$
\footnote{A \textbf{Poisson bracket} defined on the space of
real valued smooth functions on a manifold $M$, 
is a $\R$- bilinear map on $C^\infty(M)$,
$(u,v)\mapsto \{u,v\}$, such that for any $u,v,w \in C^\infty(M)$:\\
$\qquad \bullet \{u,v\}=-\{v,u\}$ (skewsymmetry),\\
 $\qquad \bullet \{u,vw\}=\{u,v\}w+\{u,w\}v $ (Leibniz rule)\\
$\qquad \bullet \{\{u,v\},w\}+\{\{v,w\},u\}+\{\{w,u\},v\}=0$ (Jacobi's identity).\\
A Poisson bracket 
is given in terms of a contravariant skew 
symmetric 2-tensor $\PO$ on $M$, called the \textbf{Poisson tensor}, by
$ \{u,v\} = \PO(du\wedge dv)$; Jacobi's identity is then equivalent to the vanisihng of the Schouten bracket $ [\PO,\PO]_S=0.$
The {\bf{Schouten bracket}} is the extension -as a graded derivation for the exterior product- of the bracket of vector fields to skewsymmetric  contravariant tensors.\\
Given a symplectic manifold $(M,\omega)$ the Poisson bracket is defined by $\{u,v\}=-\omega(X_u,X_v)$ where $X_u$ is the Hamiltonian vector field associated to $u$, i.e. $\iota(X_u)\omega=du$.
On $(\R^{2n},dp_i\wedge dq^i)$, the bracket is $\{f,g\}=\partial_{q^i}f\partial_{p_i}g-\partial_{p_i}f\partial_{q^i}g$.}
Observables are families of smooth functions on that manifold
and the dynamics is defined in terms of  a Hamiltonian $H \in C^\infty(M)$ :
 the time evolution
of an observable $\{f_t\}$ 
is governed by the equation :
$$
{\frac{d}{dt}}f_t = -\left\{H,f_t\right\}.
$$
The Heisenberg's formulation of quantum mechanics has for 
framework a Hilbert space (states are rays in that space).
 Observables are 
families   of selfadjoint operators on that Hilbert space and
the dynamics is defined in terms of a Hamiltonian $H$, which 
is a selfadjoint operator : the time evolution of an
observable $\{A_t\}$ is governed by the equation~:
$$
{\frac{dA_t}{dt}} = {\frac{i}{\hbar}} [H,A_t].
$$
A natural suggestion for quantization is a correspondence
$\mathcal{Q}\colon f\mapsto \mathcal{Q}(f)$ mapping a function $f$ to a self adjoint
operator $\mathcal{Q}(f)$ on a Hilbert space $\H$ in such a way that
$\mathcal{Q}(1)=\Id$ and
\begin{equation}\label{eq:quant}
[\mathcal{Q}(f),\mathcal{Q}(g)]=i\hbar \mathcal{Q}(\{f,g\}) +O(\hbar^2).
\end{equation}
Van Hove showed that there is no correspondence  defined
on all smooth functions on $M$ so that $$[\mathcal{Q}(f),\mathcal{Q}(g)]=i\hbar \mathcal{Q}(\{f,g\}),$$
when one puts  an irreducibility
requirement which is necessary not to violate Heisenberg's principle.
More precisely, he proved that there is 
no irreducible representation of the Heisenberg algebra, viewed as the algebra of constants and linear functions on $\R^{2n}$
endowed with  the Poisson braket, which extends to a representation of the algebra of polynomials on $\R^{2n}$.\\

A natural question is to know what would appear in the righthand side of equation \eqref{eq:quant}, i.e. what would correspond to the bracket of operators.
 Similarly,  the associative law $*$ which  would appear as corresponding to the composition of operators
\begin{equation}
\mathcal{Q}(f)\circ\mathcal{Q}(g)=\mathcal{Q}(f*g)
\end{equation}
is at the root of deformation quantization, which  expresses quantization in terms of such an associative law, without knowing a priori the map $\mathcal{Q}$. \\
In deformation quantization, the quantum observables are not constructed as usually done as operators on a Hilbert space; instead quantum and classical observables coincide; one keeps the same space of smooth functions on a Poisson manifold and quantization appears as a new associative algebra structure on this space.\\

A first step is to define such an associative law as a formal deformation of the usual product of functions giving by antisymmetrization a deformation of the Poisson bracket.
This yields the notion of a {\bf{formal deformation quantization}}, also called a (formal) star product.

\begin{definition}\cite{bib:Flato}
A  \textbf{star product} on a Poisson manifold $(M,\PO)$ is a bilinear map
\[
C^\infty (M)\times C^\infty (M) \to C^\infty (M)[[\nu]] :
 \qquad (u,v) \mapsto u\star v = u\star_\nu v : =
\sum_{r\ge0} \nu^rC_r(u,v)\]
such that :\\
(a)  when the map is extended $\nu$-linearly (and 
 continuously in the $\nu$-adic topology) to \\
 $C^\infty (M)[[\nu]]\times C^\infty (M)[[\nu]]$,
it is formally associative:
\[
(u\star v)\star w = u\star(v\star w);
\] 
(b)  $C_0(u,v) = uv=:\mu(u,v)$;\quad  $C_1(u,v)-C_1(v,u) = \{u,v\}=P(du\wedge dv)$;\\
 (c) $1\star u = u\star 1 = u$;\\[1mm]
(d) the $C_r$'s are bidifferential operators on $M$ (it is then a {\bf differential star product}).\\
 When each  $C_r$
 is of order  $\le r$ in each argument, $\star$  is called  {\bf natural }.\\
 If $\overline{ f \star  g} =\overline{ g }\star \overline{ f }$ for any purely imaginary $\nu=i\lambda$,  $\star$  is called {\bf Hermitian}.
 \end{definition}
 \begin{example} The first example is the so called  {\bf{Moyal}} or {\bf{Moyal-Weyl $\star$-product}}, defined on $V=\R^m$
endowed with a Poisson structure $ \PO=\sum_{i,j} \PO^{ij} {\partial}_{i} \wedge
 {\partial}_{j}$ with constant coefficients:
$$
 (u\star_Mv)(z)=
 \left.\exp \left(\frac{\nu}{2} \PO^{rs}
 \partial_{x^r}\partial_{y^s}\right)(u(x)v(y))\right\vert_{x=y=z}. 
$$
When $\PO$ is non degenerate, i.e. on the symplectic manifold $(V=\R^{m=2n},\omega=dp_i\wedge dq^i)$, the space of formal series of polynomials on $V$ with this Moyal  deformed product, $(S(V^*)[[\nu]],\star_M)$, 
is called {\bf the Weyl algebra}.\\
The Moyal star product  on $(V=\R^{m=2n},\omega=dp_i\wedge dq^i)$ is related  to the composition of operators via Weyl's quantisation
of  polynomials   :
\begin{eqnarray}\label{eq:weylordering}
f*_wg:&=&\mathcal{Q}_{Weyl}^{-1}\left( \mathcal{Q}_{Weyl}(f)\circ \mathcal{Q}_{Weyl}(g)\right)\cr
&=& f.g +\frac{i\hbar}{2} \{f,g\} +O(\hbar^2)=f\star_Mg\vert_{\nu={i\hbar}}
\end{eqnarray}
where the Weyl quantization $\mathcal{Q}_{Weyl}$ is the bijection  between complex-valued polynomials on $\R^{2n}$, $\C[p_i,q^j]$
and the space of differential operators with complex polynomial coefficients on $\R^n$, $D_{polyn}(\R^n)$, defined by 
 $\mathcal{Q}_{Weyl}(1)=\Id,$  $\mathcal{Q}_{Weyl}(q^i):=Q^i := q^i\cdot$ is the multiplication by $q^i$, 
$\mathcal{Q}_{Weyl}(p_i):=P_i=-i\hbar \frac{\partial}{\partial q^i}$ and to a polynomial in $p's$ and $q's$ the corresponding
totally symmetrized polynomial in $Q^i$ and $P_j$. 
 \end{example}
 
 Quantization appears in this way just as a formal deformation of a Poisson algebra $\mathcal{A}$ of classical observables. It can be formulated in this very general setting.
 The main difficulty of formal deformation quantization is that the deformation parameter $\nu$ corresponds to $i\hbar$ which is a  non zero constant, and the convergence of the formal star product has to be solved to provide a good physical model of quantization.
 Nevertheless, at the formal level, deformation quantization is a very fruitful theory.\\

 In Section \ref{section:exandclass}, we recall some of the results about existence, classification and representations  for formal star products.
We give  in Section \ref{section:groupactions} results concerning the notion of classical symmetries and invariant formal star products,  the link between  invariant formal star products on Lie groups and Drinfeld twists, and the notion of universal deformations formulas.
We then consider in Section \ref{momentmaps} actions implemented by a moment map and the notion of quantum moment in the framework of formal star products.
We describe  in Section \ref{section:reduction} one of the various quantum analogues of the classical reduction procedure which  have been considered in the framework of formal deformation quantization.  
We end up in Section \ref{section:convergence} by some considerations about convergence of star products.

 \section{Existence, classification and representations  for formal star products}\label{section:exandclass} Star products were first studied in the symplectic framework; after several classes of examples, the existence on a general symplectic manifold was proven :
 \begin{thm}[De Wilde and Lecomte, 1983 \cite{bib:DWL}]
On  any symplectic manifold  $(M,\omega)$, there exists a differential star product.
\end{thm}
 Fedosov gave in 1994 \cite{bib:F} (after a first version in Russian in 1985) a  recursive construction of such a star product on a symplectic manifold when one has chosen a 
symplectic connection\footnote{A symplectic connection is a linear connection without torsion such that the covariant derivative of $\omega$ vanishes;  such connections exist on any symplectic manifold but are not unique; indeed, given any torsion free connection ${\widetilde{\nabla}}$, one can define a symplectic connection $\nabla$ via
$\nabla_XY:={\widetilde{\nabla}}_XY+\frac{1}{3}(S(X,Y)+S(Y,X))$ with $\omega(S(X,Y),Z)=({\widetilde{\nabla}}_X\omega)(Y,Z)$. Any other symplectic connection is of the form
$\nabla'_XY=\nabla_XY+L(X,Y)$ with $\omega(L(X,Y),Z)$ totally symmetric.}  $\nabla$ and a sequence of closed $2$-forms ${\tilde{\Omega}}=\sum_{k\ge 1} \nu^k \omega_k$ on $M$.
We now briefly describe this construction; it is  obtained 
by identifying $C^\infty(M)[[\nu]]$ with the algebra of flat sections of a bundle of algebras over $M$,
 the  Weyl bundle ${\mathcal W} = F(M,\omega) \times_{Sp(V,\Omega),\rho}W$,
 endowed with a flat covariant derivative $D$ built from $\nabla$ et ${\tilde{\Omega}}$.
$F(M,\omega)$ is the bundle of symplectic frames\footnote{A symplectic frame at a point $p\in M$ is a linear symplectic isomorphism \\ $\xi: (V=\R^{2n},\Omega={\tiny{\left(\begin{matrix} 0 & I_n\\-I_n & 0\end{matrix}\right)}} )\rightarrow (T_pM,\omega_p)$; $F(M,\omega)$ is a $Sp(V,\Omega)$-principal bundle over $M$.} ; $W$ is the formal Weyl algebra which is the completion of the Weyl algebra $(S(V^*)[\nu],\star_M)$
for the grading assigning the degree $1$ to $y\in V$ and the degree $2$ to  $\nu$;  $\rho$ is  the natural representation of the symplectic group $Sp(V,\Omega)$ on $W$ extending the action on $V^*$; it acts by automorphisms of $\star_M$, which shows that  the Weyl bundle is indeed a bundle of algebras, the product in the fiber being defined by  $\star_M$.

The symplectic connection $\nabla$  induces a
covariant derivative $\partial$ of sections of  ${\mathcal W}$:\\
$ \partial a = da -
{\frac{1}{\nu}}[\half \omega_{ki} \Gamma^k_{rj} y^i y^j,a]$, with summation over repeated indices,  where the $\Gamma^k_{rj} $ are the Christoffel symbols of the connection and
where the bracket is defined on $\mathcal{W}$-valued forms by combining the
 skewsymmetrisation of  $\star_M$ on sections of $\mathcal{W}$ and exterior product of forms. 
The covariant derivative acts by derivation of the space of sections of $\mathcal{W}$, which is an algebra for the $\star_M$ pointwize product of sections with values in $W$.

One deforms the covariant derivative into $Da = \partial a - \delta(a) -
{\frac{1}{\nu}}[r,a]$ where\\
$
 \delta(a) = \frac{1}{\nu}\left[-
\omega_{ij}y^idx^j,a \right] = \sum_k dx^k\wedge
{\frac{\partial a}{\partial y^k}},$
 with $r$ a $1$-form with values in ${\cal W}$.  It is clearly still a derivation of the algebra of sections of  $\mathcal{W}$, so flat sections (i.e. sections $a$ so that $Da=0$) form a subalgebra. To have enough flat sections, one asks the covariant derivative to be flat, i.e. that its curvature $D\circ D$ vanishes.
 Now $ 
D_\circ
D a = \frac{1}{\nu}\left[\overline{R} -\partial r +\delta r +
\frac{1}{2\nu}[r,r], a \right] $ and one looks for an $r$ such that  $D_\circ
D=0$.
Such an $r$, satisfying $\delta r =-\overline{R} + \partial r -\frac{1}{\nu} r^2+{\tilde{\Omega}}$, 
can be defined inductively by $$ r =
-\hat{\delta}\overline{R} + \hat{\delta} \partial r 
- \frac{1}{\nu}
\hat{\delta} r^2+\hat\delta{\tilde{\Omega}} $$ where,
 writing  any $a \in \Gamma({\cal W}\otimes\Lambda^q)$ in the form
$$
a=\sum_{p\ge0,q\ge0}  a_{pq}= \sum_{2k+p\ge0,q\ge0} \nu^k
a_{k,i_1,\ldots,i_p,j_1,\ldots,j_q} y^{i_1}\dots
y^{i_p}\,dx^{j_1}\wedge\dots\wedge dx^{j_q},$$
one defines $\hat{\delta}(a_{pq}) =\left\{\begin{array}{l}
 \frac{1}{{p+q}}\sum_ky^ki({\frac{\partial}{\partial x^k}})a_{pq}\quad
\textrm{ if} \, \, p+q >0,\cr
0\quad
\textrm{ if} \, \,p+q=0. \end{array}\right.$\\

   A flat section of ${\mathcal W}$ is then given inductively by 
 $ a
= \hat\delta\left(\partial a -
\frac{1}{\nu} [r,a]\right) + a_{00}, $  so corresponds bijectively with $a_{00}$ which is an element of $\C^\infty(M)[[\nu]]$ and is denoted $Q(a_{00})$.
The Fedosov's star product  $*_{\nabla,\Omega}$ is then obtained  by $u\star_{\nabla,\Omega}v:=(Q(u)\star_M Q(v))_{00}$.\\
Omori, Maeda and Yoshioka gave yet another proof of existence by glueing locally defined Moyal-Weyl star products \cite{bib:OMY}.
\begin{definition}
Given a star product $\star$ and any series  $T=\sum_{r\ge 1}\nu^r T_r$ of linear operators on $\mathcal{A}=C^\infty(M)$, on can build another star product
denoted $\star':=T\bullet \star$ via \begin{equation}\label{eq:equiv}
u\star'v:=e^{T}\left(e^{-T}u \star e^{-T} v\right).
\end{equation}
Two star products $\star$ and $\star'$ are said to be {\bf{equivalent}} if there exists a series $T$ such that equation \eqref{eq:equiv} is satisfied.
If the star products are differential and equivalent, the equivalence can be  defined by a series of differential operators. 
\end{definition}
The classification of  star products up to equivalence on symplectic  manifolds was obtained  by Nest-Tsygan \cite{bib:NT},  Deligne \cite{bib:Deligne}, and Bertelson-Cahen-Gutt \cite{bib:BCG} :
 \begin{thm}
 Any star product on a symplectic manifold is equivalent to a Fedosov's one and its  equivalence class
is parametrised by the element in  $H^2(M;\R)[[\nu]]$  given by the series $[{\tilde{\Omega}}]$ of de Rham classes of the closed $2$-forms used in the construction.
 \end{thm}
Fedosov had obtained the classification of star products obtained by his construction, and Deligne gave an instrinsic way to define the characteristic class associated to a star product.
For a detailed presentation of this class, we refer to \cite{bib:GRequiv}.\\

Concerning star products on Poisson manifolds, a proof of existence  quickly followed for regular Poisson structures (by Masmoudi).  
An explicit construction of star product was known for linear Poisson structure, i.e. on the dual $\g^*$ of a Lie algebra $\g$ with the Poisson structure defined by $$\PO_\xi(X,Y):=<\xi,[X,Y]> ,\qquad \xi \in \g^*, X,Y\in \g\simeq T^*_\xi\g^*,$$using the fact that polynomials on $\g^*$ identify with the symmetric algebra $S(\g)$ which in turns is in bijection with the universal enveloping algebra  $\mathcal{U}(\g)$ which is associative.
Pulling back this associative struture to the space of polynomials on $\g^*$ yields a differential star product \cite{bib:Gutt}.\\
For general Poisson manifolds,  the problem of existence and classification of star products was solved ten years later by Kontsevich :
\begin{thm}[Kontsevich, 1995, \cite{bib:K}]\label{thm:K}
The  set of  equivalence classes of differential star products on a
 Poisson manifold $(M,\PO)$ coincides with
the set of equivalence classes of Poisson deformations of $\PO$:
 $$
 \PO_{\nu}=\PO\nu+\PO_2 \nu^2+\dots\in \nu\Gamma(X,\Lambda^2 T_X)[[\nu]],\,\,\,\textrm{such that}\,\,
 [\PO_{\nu},\PO_{\nu}]_S=0 ,
 $$ 
 where equivalence of  Poisson deformations is defined via the action of a formal vector field on $M$, $X=\sum_{r\ge 1}\nu^r X_r$, via $\{u,v\}':=e^{X}\left\{e^{-X}u , e^{-X} v\right\}$.
  \end{thm}
 Remark that in the symplectic framework, this result coincides with the previous one. Indeed any Poisson deformation $\PO_{\nu}$ of the Poisson bracket $\PO$
on a symplectic manifold $(M,\omega)$ is of the form $\PO^\Omega$
for a series $\Omega=\omega+\sum_{k\ge 1}\nu^k\omega_k$ where the $\omega_k$
are closed $2$-forms,  with
$\PO^\Omega(du,dv)=-\Omega(X^\Omega_u,X^\Omega_v),\quad X^\Omega_u
\in \Gamma (TM) [[\nu]]\, \textrm{defined by}\,
i(X^\Omega_u)\Omega=du.
$ 

We briefly sketch how Kontsevich's theorem is a consequence of his formality theorem. A general yoga sees any deformation theory  encoded in a differential graded Lie algebra structure\footnote{ 
  A {\bf differential graded Lie algebra} (briefly DGLA) is a 
graded Lie algebra $\g$ together with a differential $(\g,[\, ,\, ], d)$: 
$d \colon \g \to \g$, i.e.\ a graded derivation of degree 1 ($ d \colon \g^i \to \g^{i+1}$, $d [{a},{b}]  = [d{a},{b}] + (-1)^{\vert a \vert} [{a},{d b}]$)
so that $d\circ d =0$.\\
A {\bf{deformation}} is a Maurer-Cartan element, i.e. a   $C \in \nu\g^1[[\nu]]$ so that $dC-\half [C,C]=0.$ \\
{\bf{Equivalence}} of deformations is obtained through  the action of the group $\exp \nu {\g}^0[[\nu]]$, the infinitesimal action of a $T\in  \nu {\g}^0[[\nu]]$   being  $T\cdot C := -dC+[T,C]$.}.\\

To express star products in that framework, one considers {\bf the   DGLA of polydifferential operators}.
Let $({\mathcal{A}},\mu)$ be an associative algebra with unit on a field $\mathbb{K}$.
Consider the {\bf{Hochschild complex}} of multilinear maps from ${\mathcal{A}}$ to itself:
$
\mathcal{C }({\mathcal{A}}):= \sum_{i=-1}^{\infty} \mathcal{C }^i $ with $ \mathcal{C }^i := \Hom_{\mathbb{K}}({\mathcal{A}}^{\otimes(i +1)}, {\mathcal{A}});
$
remark that the degree is shifted by one; the degree $\vert A\vert$ of a $(p+1)$--linear map $A$ is equal to $p$.
For 
$A_1\in \mathcal{C }^{m_1}$, $A_2\in \mathcal{C }^{m_2}$, define:
\begin{eqnarray*}
&&(A_1 \circ A_2) 1(f_1,\dots, f_{m_1+m_2+1}) := \cr
&&\sum_{j=1}^{m_1}
(-1)^{(m_2)(j-1)} A_1( f_1, \ldots,f_{j-1}, A_2(f_j,\dots,f_{j+m_2}), 
f_{j+m_2+1},\ldots, f_{m_1+m_2+1}).
\end{eqnarray*}
The {\bf Gerstenhaber bracket} is defined by
$
[A_1,A_2]_G:=A_1\circ A_2-(-1)^{m_1m_2}A_2\circ A_1.
$
It gives $ \mathcal{C }$ the structure of a graded Lie algebra.
An element  $M\in \mathcal{C }^1$ defines an associative product  iff $[M,M]_G=0$.
 The {\bf{differential}} $d_\mu$ is defined by
$
d_\mu A=-[\mu,A].
$
Then $({\mathcal{C}}({\mathcal{A}}),[\, ,\, ]_G,d_\mu)$   is a differential graded Lie algebra.\\
Here  we consider ${\mathcal{A}} = C^\infty(M)$, 
and we  deal  with the subalgebra of ${\mathcal{C}}({\mathcal{A}})$
consisting of multidifferential 
operators ${\mathcal{D}}_{poly}(M):= \bigoplus {\mathcal{D}}^i_{poly}(M)$ 
with 
$ {\mathcal{D}}^i_{poly}(M)$
 the set of multi differential operators
acting on $i+1$ smooth functions on $M$
and vanishing on constants.
${\mathcal{D}}_{poly}(M)$ is closed under 
the Gerstenhaber bracket and under the differential $d_\mu$, 
so that $\left({\mathcal{D}}_{poly}(M),[\, ,\, ]_G,d_{{\mathcal{D}}}:=d_\mu\vert_{\mathcal{D}}\right)$   is a DGLA.\\[1.5mm]
A $\star$-product  is given by $\star=\mu +C$ with $C\in \nu {\mathcal{D}}^1_{poly}(M)[[\nu]]$
a Maurer-Cartan element of the DGLA $\left({\mathcal{D}}_{poly}(M),[\, ,\, ]_G,d_{\mathcal{D}}\right)$, indeed the associativity $[\mu +C,\mu +C]=0$ is equivalent to 
$
d_DC-\half [C,C]_G=0.
$\\
  Equivalence of star products is given by the  action  of $e^T$ with $T\in  \nu {\mathcal{D}}^0_{poly}(M)[[\nu]]$   via :\\  $\mu+C'=\left(\exp [ T,~]_G\right)(\mu+C)$; 
 the infinitesimal action  is $T \cdot C:= -d_{\mathcal{D}}C+[T,C]_G $.  \\

To express Poisson deformations in that framework, one considers  {\bf{the DGLA of skewsymmetric polyvectorfields}} $({\mathcal{T}}_{poly}(M),[\, ,\, ]_{\mathcal{T}},0)$
with ${\mathcal{T}}^p_{poly}(M):=\Gamma(\Lambda^{p+1} TM)$ (remark again the shift in the grading, a $p+1$-tensor field is of degree $p$); 
 the bracket is given, up to a sign, by the {\bf{Schouten bracket}} $[T_1,T_2]_{\mathcal{T}}:=-[T_2,T_1]_S$ and  the algebra is endowed with the {\bf{zero differential}} $d_{\mathcal{T}}:=0$.\\[1mm]
A Poisson deformation  is given by $\PO_{\nu}  \in \nu {\mathcal{T}}^1_{poly}(M)[[\nu]]$ so that $[\PO_{\nu},\PO_{\nu}]_S=0$, hence so that $d_{\mathcal{T}}\PO_{\nu}-\half [\PO_{\nu},\PO_{\nu}]_{\mathcal{T}}=0$ and thus by 
a Maurer-Cartan element of the DGLA.
  Equivalence of Poisson deformations  is given by the  action  of $e^T$ with $T\in  \nu {\mathcal{T}}^0_{poly}(M)[[\nu]]$   via $\PO'_\nu=\left(\exp [ T,~]_{\mathcal{T}}\right)(\PO_\nu)$; 
 the infinitesimal action  is $T \cdot \PO_{\nu}:= -d_{\mathcal{T}}\PO_{\nu}+[T,\PO_{\nu}]_{\mathcal{T}} $.  \\
 
Remark that any DGLA $(\g,[\, ,\, ], d)$ has a cohomology complex defined by
\[ 
{{H}}^i(\g) := \Ker(d \colon \g^i \to \g^{i+1}) \Big/ \Im(d \colon \g^{i-1} \to \g^i).
\]
The set ${{H}} := \bigoplus_i {{H}}^i(\g)$ 
inherits the structure of a graded Lie algebra :
$
[{|a|},{|b|}]_{{H}} := \left| [{a},{b}] \right|
$
where $|a| \in \mathcal{H}$ denote the equivalence classes of a $d$-closed element
$a \in \g$. \\
Then $(H,[\, ,\, ]_H, 0)$  is a DGLA  (with zero differential).
 \begin{thm}[Vey 1975, \cite{bib:Vey}] 
Every  cocycle $C \in {\mathcal{D}}^{p}_{poly}(M)$  (i.e. such that $d_{\mathcal{D}}(C)=0$) is the sum of the
coboundary of a $B \in {\mathcal{D}}^{p-1}_{poly}(M)$ and a $1$-differential
skewsymmetric $p$-cocycle $A$, hence
$
H^p({\mathcal{D}}_{poly}(M))=HH^p_\mathrm{diff}(C^\infty(M),C^\infty(M)) = \Gamma(\Lambda^{p+1} TM)={\mathcal{T}}^p_{poly}(M).
$
\end{thm}
The DGLA defined by the cohomology  of $\left({\mathcal{D}}_{poly}(M),[\, ,\, ]_G,d_{{\mathcal{D}}}\right)$ is $({\mathcal{T}}_{poly}(M),[\, ,\, ]_{\mathcal{T}},0)$.
The natural map $U_1 \colon  {\mathcal{T}}^i_{poly}(M)\longrightarrow {\mathcal{D}}^i_{poly}(M)$
\begin{equation}\label{eq:U1}
U_1(X_0 \wedge \ldots \wedge X_n)(f_0, \ldots , f_n)=
\frac{1}{(n+1)!}\; \sum_{\sigma \in S_{n+1}} 
\epsilon(\sigma) \; X_{{0}}(f_{\sigma(0)})\cdots X_{{n}} (f_{\sigma(n)}),
\end {equation}
intertwines the differential and induces the identity in cohomology, but is not a DGLA morphism.
A  DGLA morphism from $({\mathcal{T}}_{poly}(M),[\, ,\, ]_{\mathcal{T}},0)$ to $\left({\mathcal{D}}_{poly}(M),[\, ,\, ]_G,d_{{\mathcal{D}}}\right)$, inducing the identity in cohomology,
would give a correspondence between a formal Poisson tensor on $M$
and a formal differential star product on $M$ and a bijection between equivalence classes. 
The existence of such a morphism  fails;  to circumvent this problem, one extends
the notion of morphism between two  DGLA introducing $L_\infty$-morphisms.\\

Let $W = \oplus_{j \in \intg} W^j$ a $\mathbb{Z}$-graded vector space. Let $V=W[1]$ be the shifted graded vector space \footnote{As vector spaces, $W[1]=W$ but there is a shift in the degrees; an element of degree $p$ in $W$ has degree $p-1$ in $W[1]$.} .
The {\bf{graded symmetric bialgebra of $V$}}, denoted 
$\Sym V$, is  the quotient of the free algebra $\Tens V$ by the
two-sided  ideal generated by
$x\otimes y-(-1)^{|x||y|}y\otimes x$ for any homogeneous elements $x,y$ in $V$.
The coproduct $\Delta_{sh}$ is  induced by the morphism of associative algebras $\Delta_{sh}:\Tens V\rightarrow  \Tens V \otimes \Tens V$ so that $\Delta_{sh}(x)=1\otimes x+x\otimes 1.$\\
A {\bf{ $L_\infty$-structure on $W$}} is defined to be a graded coderivation $\mathcal{Q}$ of $\Sym\,(W[1])$ of degree $1$ satisfying $\mathcal{Q}^2=0$ and $\mathcal{Q}(\un_{\Sym W[1]})=0$. 
Such  a coderivation is determined by 
$$
Q:=pr_{W[1]}\circ \mathcal{Q}: \Sym\,(W[1])\rightarrow W[1]\, \textrm{via}\, \mathcal{Q}=\mu_{sh}\circ Q\otimes \Id \circ \Delta_{sh} \, \textrm{and we write}\,\mathcal{Q}=\overline{Q}.
$$
The pair $(W,\mathcal{Q})$ is called an {\bf{$L_\infty$-algebra}}.
\begin{example} If $(\g,[~,~],d)$ is a DGLA, then $ \left(\g, \mathcal{Q}= \overline{d[1] + [~,~][1]} \right)$\footnote{For  $\phi : V^{\otimes k}\to W^{\otimes \ell}$,  on defines  $\phi[j] :  V[j]^{\otimes k}\to W[j]^{\otimes \ell}$  via \\$\phi[j]:= (s^{\otimes \ell})^{-j}\circ \phi \circ (s^{\otimes k})^j$  where $s:V\rightarrow V[-1]$  is the suspension, i.e. the identity map with shifts of degrees. }
, with $\mathcal{Q}$ defined on  $\Sym(\g[1])$), is an $L_\infty$-algebra.
 A deformation, i.e. a $C \in \nu\g^1[[\nu]]$ so that  $dC+\half[C,C]=0$, corresponds to a series $C\in \nu(\g[1])^0[[\nu]]$ such that $\mathcal{Q}(e^{C})=0$.
\end{example}
A {\bf{$L_\infty$-morphism}} from a $L_\infty$-algebra 
$(W,\mathcal{Q})$ to a $L_\infty$-algebra $(W',\mathcal{Q}')$ is a morphism of graded 
connected  coalgebras $\Phi :\Sym(W[1]) \to \Sym(W'[1])$, intertwining differentials  $ \Phi \circ \mathcal{Q} = \mathcal{Q}' \circ \Phi$.
Such  a morphism is determined by $$\varphi :=pr_{W'[1]}\circ \Phi: \Sym\,(W[1])\rightarrow W'[1]\,\textrm{ with }\,\varphi(1)=0,\,\textrm{ via }\, \Phi=e^{*\varphi}$$
with $A*B=\mu \circ A\otimes B\circ \Delta$ for $A,B\in \Hom(\Sym(W[1]),\Sym(W'[1]))$.\\
$\Phi$ is  a {\bf{quasi-isomorphism}} if  $\Phi_1= \Phi|_{W[1]}=\varphi_1:W[1]\to W'[1]$  induces an isomorphism  in cohomology. \\

A {\bf{formality for a DGLA}} $(\g,[~,~],d)$ is a
quasi-isomorphism 
from the $L_\infty$-algebra corresponding to $(H,[~,~]_H,0)$ (the cohomology of $\g$ with respect to $d$ ), 
to the $L_\infty$-algebra corresponding to $(\g,[~,~],d)$; it is thus a map
$$\Phi: \Sym (\mathfrak{H}[1])\rightarrow \Sym (\g [1]) \,{\textrm{ such that }}\, \Phi \circ \overline{[~,~]_H[1]}=\overline{(d[1]+[~,~][1])}\circ \Phi.$$
In case one has a formality, the space of deformations modulo equivalence coincide for  $(\g,[~,~],d)$ and $(H,[~,~]_H,0)$.
In particular, a formality  for $\left({\mathcal{D}}_{poly}(M),[\, ,\, ]_G,d_{{\mathcal{D}}}\right)$ yields a proof of theorem \ref{thm:K}.
For a Poisson structure $\PO$, the associated Kontsevich star product   is given by $\star_{\PO}=\mu +\sum_{k\ge 1} \nu^{k}F_k(X,\PO, \cdots ,\PO)$
where the $F_k$ are the so-called Taylor coefficients, i.e. projections  on ${\mathcal{D}}_{poly}(\R^d)[1]$ of the formality restricted to $S^k({\mathcal{T}}_{poly}(\R^d)[1])$.\\
 
Kontsevich gave an explicit formula for  a formality for $\left({\mathcal{D}}_{poly}(M),[\, ,\, ]_G,d_{{\mathcal{D}}}\right)$ when $M=\R^d$: he gave the Taylor coefficients $F_n$ of an 
$L_\infty$--morphism  between the two $L_\infty$-algebras 
$$
F: ( {\mathcal{T}}_{poly}(\R^d),\mathcal{Q})\rightarrow ({\mathcal{D}}_{poly}(\R^d),\mathcal{Q}')
$$
corresponding  to the  DGLA's  $ ({\mathcal{T}}_{poly}(\R^d)~,~[~Ê,~]_{\mathcal{T}},~d_{\mathcal{T}}=0)$
and   $({\mathcal{D}}_{poly}(\R^d)~,~[~,~]_G~,~d_{\mathcal{D}})$
with the first coefficient 
$
F_1: {\mathcal{T}}_{poly}(\R^d)\rightarrow {\mathcal{D}}_{poly}(\R^d)
$
given  by $F_1=U_1$ as in \eqref{eq:U1}.
The formula is of the form
\[
F_n=\sum_{m\ge 0}\sum_{{\vec{\Gamma}}\in G_{n,m}}{\mathcal{W}}_{\vec{\Gamma}}B_{\vec{\Gamma}}
\]
where $G_{n,m}$ is a set of oriented admissible graphs;
$B_{\vec{\Gamma}}$ associates a $m$--differential operator
to an $n$--tuple of multivectorfields;
and ${\mathcal{W}}_{\vec{\Gamma}}$ is the integral 
of a form $\omega_{\vec{\Gamma}}$ over the compactification
of a configuration space $C^+_{\{p_1,\ldots,p_n\}\{q_1,\ldots,q_m\}}$. For details, we refer to \cite{bib:K,bib:Arnal}\\
An explicit globalisation on a manifold has been built by Cattaneo, Felder and Tomassini \cite{bib:CFT}, who also gave an interpretation of the formula in terms of sigma models \cite{bib:CF}.\\
Given a manifold and a torsion free connection $\nabla$ on it, Dolgushev \cite{bib:Dolgushev} has built a formality
$
F: ( {\mathcal{T}}_{poly}(M),\mathcal{Q})\rightarrow ({\mathcal{D}}_{poly}(M),\mathcal{Q}').
$\\

\smallskip

An interesting feature of formal deformation quantization is the possibility to define a notion of states and to study representations of the deformed algebras. For this, parts of the algebraic theory of states and representations which exist for $C^*$-algebras\footnote{ A $C^*$-algebra is a Banach algebra over $\C$ endowed with a $*$ involution (i.e. an involutive 
semilinear antiautomorphism) such that  $\Vert a \Vert = \Vert a^*\Vert$ and $\Vert a  a^*\Vert=\Vert a \Vert^2$ for each element $a$ in the algebra.
If $\mathcal{A} = \mathcal{B}(\mathcal{H})$ is the algebra of bounded linear operators on a Hilbert space 
$\mathcal{H}$ and if $\psi$ is a non vanishing element of  $\mathcal{H}$, the ray it generates
defines the linear functional $$\omega: \mathcal{A}\rightarrow \C: 
A\mapsto \omega(A):=\frac{<\psi,A\psi>}{<\psi,\psi>}$$
which is positive in the sense that $\omega(A^*A)\ge 0$.
This lead to define a state in the theory of $C^*$ algebras as a positive linear functional.} 
have been extended by Bordemann, Bursztyn  and Waldmann  \cite{bib:BW, bib:BurzsWald}   to the framework of $*$-algebras over ordered rings\footnote{.
An associative commutative unital ring $R$ is said to be {\bf{ordered}} with positive elements $P$ if
the product and the sum of two elements in $P$ are in $P$,
and if $R$ is the disjoint union $R=P\cup \{0\}\cup -P$. Examples are given by $\Z,\mathbb{Q},\R,\R[[\lambda]]$; in the case of $\R[[\lambda]]$, a series $a=\sum_{r=r_0}^\infty a_r\lambda^r$ is positive
if its lowest order non vanishing term is positive $ (a_{r_0}>0)$.}.\\
 Let $R$ be an ordered ring and $C=R(i)$  be the ring extension by a square root $i$ of $-1$ (for deformation quantization, $C=\C[[\lambda]]$
for  $R=\R[[\lambda]]$ with $\nu=i\lambda$). \\
An associative algebra  $\mathcal{A}$ over $C$  is called a {\bf{ $~^*$-algebra}} if it has an involutive antilinear antiautomorphism
$~^*: \mathcal{A} \rightarrow\mathcal{A}$ called the $~^*$-involution; for instance  $(C^\infty(M)[[\nu=i\lambda]],\star)$ with a Hermitian star product and conjugaison is a $~^*$-algebra over $\C[[\lambda]]$.\\
 A linear functional $\omega:\mathcal{A}\rightarrow C$ over a $~^*$-algebra over $C$ is called {\bf{positive}} if 
\begin{equation*}
\omega(A^*A)\ge 0 \qquad \textrm{ for any } A \in
\mathcal{A}.
\end{equation*}
 A {\bf{state}} for a $~^*$-algebra $\mathcal{A}$ with unit  over $C$ is a positive linear functional so that $\omega(1)=1$.\\
The positive linear functionals on $C^\infty(M)$ are the compactly supported Borel measures.
The $\delta$-functional on $\R^{2n}$ is not positive with respect to the Moyal star product : if $H:=\frac{1}{2m}p^2+kq^2$,
$\left({H}\star_{M}H\right)(0,0)=\frac{k\nu^2}{2m}=\frac{-k\lambda^2}{2m}<0.$
Bursztyn and Waldmann  proved in \cite{bib:BurzsWald} that for a Hermitian star product, any classical state $\omega_0$ on $C^\infty(M)$ can be deformed into
a state for the deformed algebra, $\omega=\sum_{r=0}^\infty \lambda^r\omega_r$.\\
Given a positive functional $\omega$ over the $~^*$-algebra ${\mathcal{A}}$, one can extend  the GNS construction of an associated representation of the algebra: the {\bf{Gel'fand ideal of $\omega$}} is ${\mathcal{J}}_\omega = \left\{a
    \in {\mathcal{A}} \; \big| \; \omega(a^*a) = 0 \right\}$ and on obtains the {\bf{GNS- representation}} of the algebra 
${\mathcal{A}}$ by left multiplication on the space $\mathcal{H}_\omega = { \mathcal{A}}
\big/ { \mathcal{J}}_\omega$ with the pre Hilbert space
structure defined via $\SP{[a], [b]} = \omega(a^*b)$ where $[a]=a+{ \mathcal{J}}_\omega$ denotes the class in ${ \mathcal{A}}
\big/ { \mathcal{J}}_\omega$ of $a\in { \mathcal{A}}$.\\
In that setting, Bursztyn and Waldmann introduced a notion of strong Morita equivalence (yielding equivalence of $*$-representations) and the complete classification of star products up to Morita equivalence was given, first on a symplectic and later in collaboration with Dolgushev on a general Poisson manifold \cite{bib:BDW}.\\

Another success in formal deformation quantization is the algebraic index theorem (which will not be presented here). It is an adaptation of the algebraic part of the Atiyah-Singer index theorem from pseudo-differential operators to a more general class of deformation quantizations (observing that  pseudo-differential operators on a manifold $N$ form, via their symbols and composition,  a deformation quantization of the cotangent bundle $T^*N$), the algebraic input entering the index theorem being  the equality of certain cyclic cocycles.    This algebraic index theorem was obtained by  Fedosov \cite{bib:Fed}, Nest and Tsygan \cite{bib:NT} on symplectic manifolds
and Dolgushev and Rubtsov \cite{bib:DR} for Poisson manifolds.

 \section{Symmetries and invariant formal star products}\label{section:groupactions}
 
In the framework of classical mechanics, symmetries appear in the following way.
A Lie group $G$ is a symmetry group for our classical system $(M,\PO)$ if  it acts by {\bf{Poisson diffeomorphisms}} on  $(M,P)$, i.e. iff
$$\{\,  g^*u,g^*v\,\}=g^*(\{u,v\})$$ for all $ u,v \in C^\infty(M)$ and $ g\in G,$
 or, equivalently, if and only if $g_*P=P$ for all $g\in G$.
 In the symplectic case, this is equivalent to $g^*\omega=\omega$ for all $g\in G$.\\
Any $X$ in the Lie algebra $\g$ of $G$
gives rise to a fundamental vector field $X^{*M}$ defined by
$
X^{*M}_p=\frac{d}{dt}_{\vert 0} \exp -tX \cdot p;
$
signs have been chosen so that
$
[X^{*M},Y^{*M}]=[X,Y]^{*M}$ so that one has a morphism from $\g$ into the space of vector fields $\chi(M)$ . One has an infinitesimal {\bf{Poisson action of the Lie algebra $\g$}:
\begin{equation}
\mathcal{L}_{X^{*M}}\{u,v\}=\{\mathcal{L}_{X^{*M}}u,v\}+\{u,\mathcal{L}_{X^{*M}}v\}
\end{equation}}
or equivalently
$\mathcal{L}_{X^{*M}} P=0$. In the symplectic case, this is  equivalent to $\mathcal{L}_{X^{*M}}\omega=0$ which says that
$\iota({X^{*M}})\omega$ is a closed $1$-form on $M$ for any element $X\in\g$.\\

The action of a Lie group on the classical Hilbert space framework of quantum mechanics is described by a
unitary representation of the group on the Hilbert space; such a representation acts by conjugaison on the set of selfadjoint operators on that space and yields an automorphism of the algebra of quantum observables. \\

To define symmetries in the setting of deformation quantization,  one first observes that the classical action of a group $G$ on a Poisson manifold
extends by pullbacks to an action of the space of functions and thus to the algebra of observables $C^\infty(M)[[\nu]]$ and one can define  different notions
of invariance of the deformation quantization under the action of a Lie group.

Let $(M,P)$ be a Poisson manifold, $G$ be a Lie group acting on $M$, and 
$(\C^\infty(M)[[\nu]],\star)$ be a deformation quantization of $(M,P)$.
The star product is said to be {\bf {geometrically invariant}} if, 
$$
g^*\left( u\star v\right)= g^*u \star g^*v  \qquad \forall g\in G, \forall u,v\in C^\infty(M).
$$
This clearly implies  that $
g^*\left(\{ u,v\}\right)=\{ g^*u,g^*v\}
$
so $G$ acts by  Poisson diffeomorphisms. Each fundamental vector field  $X^{*M}$ is then a derivation of the star product
$$
\mathcal{L}_{X^{*M}}\left( u\star v\right)= (\mathcal{L}_{X^{*M}}u) \star v + u \star (\mathcal{L}_{X^{*M}}v).
$$
Symmetries in quantum theories are automorphisms of the algebra
of observables; a symmetry $\sigma$ of a star product
$\star$ is an automorphism of the $\C\bbnu$-algebra $(C^\infty(M)\bbnu,\star)$ 
\[
\sigma(u\star v) = \sigma(u) \star  \sigma(v), \qquad \sigma(1) =1,
\]
where $\sigma$ is  a formal series of linear maps. One can show that  $\sigma(u) = T(u \circ
\tau)$ where $\tau$ is a Poisson diffeomorphism of $(M,P)$ and $T = \Id
+ \sum_{r\ge1} \nu^r T_r$ a formal series of differential maps. 
A Lie group $G$ acts as  {\bf {symmetries of }} $(C^\infty(M)[[\nu]],\star)$ if there
is a homomorphism $\sigma: G\rightarrow \Aut (M,\star).$ In that case, 
$\sigma(g) u=T(g)( \tau(g)^*u)$  and
$\tau$ defines a Poisson action of $G$ on $(M,P)$.\\

The existence and classification of invariant star products on a Poisson manifold is known, provided there exists an invariant connection on the manifold.
To define the equivalence in this context,
 two $G$-invariant star products  are called {\bf{$G$-equivalent} }if there is an equivalence
 between them which commutes with the action of $G$.

 Fedosov's construction in the symplectic case   builds a star product $\star_{\nabla,0}$, knowing a symplectic connection $\nabla$.
 If that connection is invariant under the action of $G$ it is clear that the construction yields an invariant star product.
 More generally,  any diffeomorphism $\phi$ of $(M,\omega)$  is a symmetry of the Fedosov star product $\star_{\nabla,{\tilde{\Omega}} }$ iff it preserves
the symplectic $2$-form $\omega$, the connection $\nabla$ and the series of closed $2$-forms $\tilde\Omega$.

 Reciprocally, we showed \cite{bib:GR}  that a natural star product on a symplectic manifold determines in a unique way a symplectic connection.
 Hence, when $G$ acts on  $(M,\omega)$ and leaves a  natural $\star$ product invariant,  there is a unique symplectic connection 
 which is invariant under $G$.
 
\begin{thm}[Bertelson, Bieliavsky, G. \cite{bib:BBG}]
Suppose $\star$ is $G$-invariant on $(M,\omega)$ and assume there exists a $G$-invariant symplectic connection $\nabla$.
 Then, there exists a series of $G$-invariant closed $2$-form $\Omega\in Z^2(M;\R)^{G-inv}[[\nu]]$ such that $\star$ is $G$-equivalent
 to the Fedosov star product constructed from  $\nabla$ and $\Omega$.
 Furthermore $*_{\nabla,\Omega}$ and $*_{\nabla,\Omega'}$ are $G$- equivalent if and only if
 $\Omega-\Omega'$ is the boundary of a series of $G$-invariant $1$-forms on $M$.
 Hence there is a bijection between the $G$-equivalence classes of $G$-invariant  $*$-products on $(M,\omega)$
 and the space of formal series of elements in the second space of invariant cohomology of $M$, $H^2(M,\R)^{G-inv}[[\nu]].$
\end{thm}
 Using Dolgushev's construction \cite{bib:Dolgushev} of a formality starting from a connection, one has a similar result in the Poisson setting : 
 \begin{thm}[Dolgushev]
 If there exists an invariant connection, there is a bijection between the $G$-equivalence classes of $G$-invariant  $*$-products on $(M,\PO)$  and the $G$-equivariant equivalence classes of $G$-invariant  Poisson deformations of $\PO$.
 \end{thm}
 Let us mention that there exist symplectic manifolds which are $G$-homogeneous but do not admit any $G$-invariant symplectic connection. A first example was given by Arnal: the orbit of a filiform nilpotent Lie group in the dual of its algebra.\\
 
The class of manifolds with a simply transitive  action  are Lie groups with the action given by left multiplication;  one is interested in left invariant $\star$-products on  Lie groups. 
Since left invariant differential operators on a Lie group $G$ are identified with elements in the universal enveloping algebra $\mathcal{U}(\g)$,
bidifferential operators can be viewed as elements of $ \mathcal{U}(\g)\otimes \mathcal{U}(\g)$ and
a left invariant $\star$-product on a Lie group $G$ is given by an element $F\in \left( \mathcal{U}(\g)\otimes \mathcal{U}(\g)\right)[[\nu ]]$,
such that
\begin{itemize}
\item $( \Delta\otimes \id)(F)\circ (F\otimes 1)=(\id\otimes \Delta)(F)\circ (1\otimes F)$ where $\circ$ denotes the product in \\$ \mathcal{U}(\g)\otimes \mathcal{U}(\g)\otimes \mathcal{U}(\g)$ and $\Delta:  \mathcal{U}(\g)\rightarrow  \mathcal{U}(\g)\otimes \mathcal{U}(\g)$ is the usual coproduct\footnote{$\Delta:  \mathcal{U}(\g)\rightarrow  \mathcal{U}(\g)\otimes \mathcal{U}(\g)$ is the algebra morphism such that $\Delta(x)=1\otimes x+ x\otimes 1$ for $x\in \g$.} , both extended $\C[[\nu]]$-linearly ; this expresses the associativity;
\item $(\epsilon\otimes \id)F=1=(\id \otimes \,  \epsilon) F$, where $\epsilon :\mathcal{U}(\g)\rightarrow \C $ is the counit; this expresses that $1\star u=u\star 1=u$;
\item $ F=1\otimes 1+ O(\nu)$, which expresses that the zeroth order term is the usual product of functions.
\end{itemize}
Such an element is called a {\bf{formal Drinfeld twist}}.
The skewsymmetric part of the first order term, which is automatically in $\g\otimes \g$ corresponds to a left invariant Poisson structure on $G$ and is what is called a {\bf{classical $r$-matrix.}}
An invariant equivalence is given by an element $S \in  \mathcal{U}(\g)[[\nu ]]$ of the form $S=1+O(\nu)$ and the equivalent $\star$-product is defined by the new Drinfeld twist  given by
$$F'=\Delta(S^{-1})F(S\otimes S).$$
Drinfeld has proven in 83 that any classical $r$-matrix arises as the first term of a Drinfeld twist (see also Halbout  about formality of bialgebras \cite{bib:H},
or  Esposito, Schnitzer and Waldmann in 2017 about a universal construction \cite{bib:ESW} ).
An  analogous algebraic construction on a homogeneous space $M=G/H$ was given by Alekseev and Lachowska in 2005 \cite{bib:AL}; invariant bidifferential operators on $G/H$ are viewed as elements of  $\left(({\mathcal{U}}(\g)/{\mathcal{U}}(\g)\cdot \h)^{\otimes2}\right)^H$; a star product is given in terms of a series  $B \in \left(({\mathcal{U}}(\g)/{\mathcal{U}}(\g)\cdot \h)^{\otimes2}\right)^H[[\nu]]$ and associativity writes again as $((\Delta\otimes \id)B) (B\otimes 1)=((\id\otimes\Delta )B)(1\otimes B)$ where both sides define uniquely invariant tri-differential operators on $G/H$.\\


Given a left invariant star product on a Lie group, hence a formal Drinfeld twist $F\in \left( \mathcal{U}(\g)\otimes \mathcal{U}(\g)\right)[[\nu ]]$ on its Lie algebra $\g$, one can deform any associative algebra $(A,\mu_A)$ acted upon by $\g$ through derivations.
This process is called a {\bf{universal deformation formula}} and is defined as follows:
$$
a\star_F b:= \mu_A\left( F\bullet (a\otimes b)\right)
$$ 
where $\bullet$ denotes the action of  $\mathcal{U}(\g)\times  \mathcal{U}(\g)[[\nu]]$ on $A\times A[[\nu]]$ which is  the  extension of the action of $\g$ on $A$ to an action of $\mathcal{U}(\g)\times  \mathcal{U}(\g) $ on $A\times A$ extended $\C[[\nu]]$-linearly.
The properties of a twist ensure that $\star_F$ is an associative deformation of $\mu_A$.
Those were studied in particular  by Giaquinto and Zhang in  \cite{bib:GZ},  by Bieliavsky and Gayral in  \cite{bib:BiGa}  in a non formal setting, and by Esposito et al in \cite{bib:ESW}.\\

An equivariant version of the algebraic index theorem, considering the action of a discrete group on a formal deformation quantization on a symplectic manifold, was obtainned by Gorokhovsky, de Kleijn and Nest in \cite{bib:GKN}.

\section{ Classical and quantum Quantum moment maps}\label{momentmaps}
 
 Of particular importance in physics is the case where the action is implemented by a moment map.
Recall that an action of a Lie group is called {\bf{(almost) Hamiltonian}} when each fundamental vector field is Hamiltonian, i.e. when 
for each $X\in\g$ there exists a function $f_X$ on $M$ such that
$$
X^{*M}u=\{f_X,u\} \qquad \forall u\in C^\infty(M).
$$
In the symplectic case this amounts to say that $\iota({X^{*M}})\omega=df_X$.
When the Hamiltonian governing the dynamics on $(M,P)$ is invariant under the action of $G$, any of those functions
$f_X$ is a constant of the motion.
A further assumption is to ask that the fundamental vector fields are Hamiltonian by means of an $G$-equivariant map
from $M$ into the dual of the Lie algebra, $G$ acting on $\g^*$ by $Ad^*$,  $J: M\rightarrow \g^*:p\mapsto J(p)$, i.e.
$$
 X^{*M}u=\{J_X,u\} \quad \forall u\in C^\infty(M) \quad \textrm{ with} \quad J_X(p)=<J(p),X>
$$
 where $<.,.>$ denotes the pairing between $\g$ and its dual.
One says then that  the action possesses {\bf{a $G$ equivariant moment map}} $J$.
Equivariance means that the Hamiltonian functions $J_X$ satisfy $J_X(g\cdot p)=J_{Adg^{-1}X}(p)$ and thus
$$
\{\,J_X,J_Y\,\}=J_{[X,Y]} \qquad \forall X,Y\in \g.
$$
An action so that each fundamental vector field is Hamiltonian and so that the correspondence $X\mapsto f_X$ can be 
chosen to be a homomorphism of Lie algebras is also called a {{\bf{strongly Hamiltonian action}}. When the group $G$ is connected, it is equivalent to the existence of a $G$
equivariant moment map, with $J_X=f_X$.

In the framework of deformation quantization, this translates in the following notions.
{\bf{An action of the Lie algebra $\g$ on the deformed algebra}}, $\left( C^\infty(M)[\nu]]],\star \right)$,
is a homomorphism $D :\g\rightarrow \Der( M,\star)$ into the space of derivations of the star product. 
A derivation $D$  is  \textbf{essentially inner}
or \textbf{Hamiltonian} if $D = \frac1\nu \ad_\star u$ for some $ u \in
\CiM\bbnu$.
We call an action of a Lie algebra (or of a Lie group) on a deformed algebra
\textbf{almost $\star$-Hamiltonian} if each $D(X)$, for any $X\in \g$, is essentially inner, and
we call (quantum) Hamiltonian a linear choice of functions $u_X \in \CiM\bbnu$  satisfying
\[
D(X) = \textstyle{\frac1\nu} \ad_\star u_X:=\textstyle{\frac1\nu} (u_X\star \cdot -\cdot \star u_X), \qquad X \in \g.
\] 
The action is \textbf{$\star$-Hamiltonian} if the formal functions $u_X$ can be chosen to
make the map 
$$ \g \to C^\infty(M)\bbnu : X \mapsto u_X$$
 a homomorphism of Lie
algebras. \\
When $\star$ is invariant under the action of $G$ on $(M,\PO)$ and  the corresponding action
of the Lie algebra $\g$ given by $D(X)=X^{*M}$ is $\star$-Hamiltonian,   
a map $ \g \to C^\infty(M)\bbnu$ as above is called \textbf{a quantum moment
map}. It is  a homomorphism of algebras such that
\[
X^{*M} = \textstyle{\frac1\nu} \ad_\star u_X, \qquad X \in \g.
\]  
For a strongly Hamiltonian action  of a Lie group $G$ on $(M,P)$,
a star product is said to be {\bf {covariant}} under $G$ if
$
f_X \star f_Y-f_Y\star f_X=\nu f_{[X,Y]} \, \forall X,Y\in\g
$
where $f:\g\rightarrow \C^\infty(M)$ is the homomorphism describing the classical moment map (i.e. $X^{*M}u=\{f_X,u\}$) and 
a star product  is called {\bf{strongly invariant}} if it is geometrically invariant and if $$
f_X \star u-u\star f_X=\nu \{ f_X,u\} \, \forall X\in\g, \forall u \in C^\infty(M).$$
 (Observe that the second condition implies that the star product is invariant under the action of the connected component of the identity in $G$.) In that case, $f$ is a quantum moment map.\\

\begin{prop}[G. -Rawnsley \cite{bib:GR}, Bahns-Neumaier \cite{bib:BN}, Kravchenko \cite{bib:Kra}]
A vector field $X$ on $M$  is a derivation of the Fedosov star product $\star_{\nabla,\tilde\Omega}$ iff
$\L_X \omega=0$, $\L_X \tilde\Omega=0$, and $\L_X \nabla=0$.
 This vector field $X$ is an inner derivation of $\star_{\nabla,\Omega}$
iff  $\L_X \nabla=0$ and there exists a series of functions $\lambda_X\in C^\infty(M)[[\nu]]$ such that
$
i(X)\omega-i(X)\tilde\Omega = d \lambda_X.
$
In this case $$
X(u)=\textstyle{\frac1\nu}(\ad_\star \lambda_X )(u).$$

A $\g$ -invariant Fedosov star product for $(M, \omega)$ is obtained from 
 a $\g$ invariant connexion and a $\g$ invariant series of closed $2$-forms $\Omega$. It
admits a quantum Hamiltonian if and only if there is a linear map $ J : \g\rightarrow C^\infty(M)[[\nu]]$ such that
 $$
dJ(X) = \iota(X^{*M})\omega - \iota(X^{*M})\tilde\Omega   \qquad \forall X\in \g.
 $$
We then have $X^{*M}u=\frac{1}{\nu}\ad_\star J(X) u.$ 
 It admits a quantum moment map if, furthermore, the linear map $ J : \g\rightarrow C^\infty(M)[[\nu]]$  can be chosen so that 
  $$
   J([X,Y])=-\omega(X^{*M},Y^{*M}) +\Omega(X^{*M},Y^{*M})\qquad \forall   X,Y \in \g.
   $$
Any symplectic manifold $(M,\omega)$ equipped with a $\g$-strongly hamiltonian action with moment map $J$ and a $\g$-invariant connection, admits 
   strongly invariant star products.
\end{prop}
 
  If one considers a pair $(\star,J)$ of an $\g$-invariant star-product and a quantum moment  map, there is a natural notion of equivalence : two such pairs $(\star,J)$ and $(\star',J')$, are ``equivariantly'' equivalent if there is a $\g$-invariant equivalence $T$ between $\star$ and $\star'$ such that $J'=TJ$. The following result gives a link with equivariant cohomology\footnote{
  Let $M$ be a manifold, $\g$ be a Lia algebra, and $\rho: \g \to \chi(M)$ be a Lie algebra morphism (i.e. an action of $\g$ on $M$). The complex of $\g$-equivariant forms $\Omega_\g(M)$ is defined as
  $$
  \Omega_{\g} (M):=\left(\oplus_{2i+j=k} [ S^i (\g^*)\otimes {\Omega}^j (M) ]^{\g}   \,  ,\ , d_{\g}:= d + \iota_{\bullet} \right)
  $$
  where invariants are taken with respect to $X\cdot ( P \otimes \beta)=(\ad^*_X P) \otimes \beta + P \otimes (\L_{\rho(X)} \alpha)$ and 
  where 
  $( \iota_{\bullet}( \alpha)) (X) =  \iota (\rho(X))(\alpha(X))$ when $ \alpha$ is viewed as a polynomial on $\g$ with values in  $\Omega(M)$.
  }.
  \begin{prop}[Reichert-Waldmann, 2017 \cite{bib:RW}]
On any symplectic manifold $(M,\omega)$  with a $\g$-strongly hamiltonian action with moment map $J$, admitting a $\g$-invariant connection, the ``equivariant" equivalence classes of pairs 
$(\star,J)$ (of an $\g$-invariant star-product and a quantum moment  map) are parametrized by series of second equivariant cohomology classes ($\frac{\vert \omega -J \vert}{\nu}+H^2_\g(M)[[ \nu ]]$) 
  \end{prop}
 
Concerning   a Kontsevich star product   $\star_{\PO}=\mu +\sum_{k\ge 1} \nu^{k}F_k(X,\PO, \cdots ,\PO)$,
defined for a Poisson structure $\PO$, given a vector field $X$ so that $\L_X\PO=0$, then   
$$A_X=X+\sum_{k\ge 1} \nu^{k}F_{k+1}(X,\PO, \cdots ,\PO)
$$
is automatically a derivation of  $\star_{\PO}$.
If $X,Y$ are two vector fields  on $M$ preserving $\PO$ then 
$$
[A_{X},A_{Y}]=A_{[X,Y] } + \sum_{k\ge 1} \nu^{k}F_{k+2} (X,Y,\PO, \cdots ,\PO).
$$
Recently, Esposito, de Kleijn and Schnitzer have proven in \cite{bib:EKS} an equivariant version of formality of multidifferential operators for a proper Lie group action; this allows to obtain a quantum moment map from a classical moment map with respect to a $G$-invariant Poisson structure  and generalizes the theorem cited above from the symplectic setting to the Poisson setting.\\

A natural class of symplectic manifolds on which there is a strongly hamiltonian action of a Lie group is the class of coadjoint orbits of Lie groups in the dual of their Lie algebras.
Much work has been devoted to the construction of interesting star-products on these orbits. The star product defined on $\g^*$ does in general not restrict to the orbits of $G$ in $\g^*$.
In fact those orbits do not always possess an invariant connection so one can not hope to get in all cases an invariant star-product. \\
For a nilpotent Lie group, Arnal and Cortet \cite{bib:AC} have  built a covariant star product using Moyal star product in good adapted coordinates. They showed that a covariant star product gives rise to a representation of the group into the automorphisms of the star product. One can define the star exponential of the elements in the Lie algebras, and this gives a construction of adapted Fourier transforms \cite{bib:AG}. They extended their construction to orbits of exponential solvable groups. \\
On the orbits a compact group $G$,   a formal star product was obtained in \cite{bib:CGR} 
by an asymptotic expansion of the associative product given by the translation at the level of Berezin's symbols of the composition of operators naturally defined by geometric quantization  in the finite  Hilbert spaces of sections of powers of a line bundle built on the K\" ahler manifold $G/T$. We recall here this notion of  Berezin's symbol.\\
Let $(L\stackrel{\pi}{\rightarrow} M,\nabla,h)$ be a quantization bundle over the compact K\"ahler  manifold $(M,\omega,J)$ (i.e.,
 $L$ is a holomorphic line bundle with connection $\nabla$ admitting an
invariant hermitian structure $h$, such that the curvature is curv$(\nabla)
= -2i\pi\omega$). Let $\H$ be the Hilbert space of holomorphic sections of
$L$.\\
 Since evaluation at a point is a continuous linear functional on $\H$, let, for any $q\in L_0$ be $e_q$ the so-called {\bf{ coherent state}}  defined by $$s(\pi(q))=<s,e_q>q \quad \textrm{for any } \, s\in \H;$$ then $e_{cq}=\overline{c}^{-1}e_q$ for any $0\neq c\in \C$, and  let $\epsilon$ be the {\bf{ characteristic function}}  on $M$ defined by $\epsilon(x)=\Vert q \Vert^2\Vert e_q \Vert^2$, with $q \in L_0$ so that $\pi(q)=x$.\\ 
 Any  linear operator ${\bf{A}}$ on $\H$ has a {\bf{Berezin's symbol}} 
\begin{equation}\label{eq:Bsymbol}\hat A(x):=\frac{<{\bf{A}}e_q,e_q>}{\Vert e_q \Vert^2}\quad  q \in L_0, \pi(q)=x\in M
\end{equation}
which is a real analytic function on $M$. The operator can be recovered from its symbol:
$$({\bf{A}}s)(x):=\int_M h_y(s(y),e_q(y)){\hat{A}}(x,y)\frac{\omega^n(y)}{n!} q\qquad s\in \H, q,q' \in L_0, \pi(q)=x, \pi(q')=y,$$ 
where $ {\hat{A}}(x,y):=\frac{<{\bf{A}}e_{q'},e_q>}{<e_{q'},e_q>}$ is the analytic continuation of the symbol, holomophic in $x$ and antiholomorphic in $y$, defined on the open dense set of $M\times M$
consisting of points $(x,y)$ such that $<e_{q'},e_q>\neq 0$.
Denote by ${\hat E}(L)$ the space
of these symbols.\\
 For any positive integer $k$,
$(L^k=\otimes^kL,\nabla^{(k)}, h^{(k)})$ is a quantization bundle for $(M,
k\omega, J)$. If ${\H}^k$ is the Hilbert space of holomorphic sections of
$L^k$, we denote by $\ehat[k]$ the space of symbols of linear operators on
${\H}^k$. If, for every $k$, the characteristic function
$\epsilon^{(k)}$ on $M$
 is constant (which is true in a homogeneous case), one says that the quantization is {\bf{ regular}}. In that case, the space $\ehat[l]$ is contained in the space
$\ehat[k]$ for any $k\ge l$. Furthermore  $\CL:=\cup_{l=1}^\infty \ehat[l]$
 is a dense subspace of the space of continuous functions
on $M$. Any function $f$ in $\CL$ belongs to a particular $\ehat[l]$ and is
thus the symbol of an operator $ {\bf{A}}_f^{(k)}$ acting on ${\H}^k$ for $k\ge
l$. One has thus constructed, for a given $f$, a family of quantum
operators parametrized by an integer $k$.
From the point of view of deformation theory,  one has constructed a family of associative products $*_k$ on
$\ehat[l]$, with values in $\CL$, parametrized by an integer $k$ :
\begin{equation} \label{comp}
f *_k g=\widehat{{\bf{A}}_f^{(k)}{\bf{A}}_g^{(k)}} \qquad f, g \in \ehat[l];\,k\ge
l.
\end{equation} Similar methods were developed in the framework of Toeplitz quantization by Bordemann, Meinrenken, Schlichenmaier and Karabegov \cite{bib:BMS,bib:KS} , including operator norm estimates to obtain a continuous field
of $C^*$-algebras.\\
We  gave an algebraic construction of a star product on polynomials restricted to some orbits of a semisimple Lie group, but those star products are usually not differential.\\
 
 \section{Marsden-Weinstein reduction of a strongly invariant star product}\label{section:reduction}
 Reduction is an important classical tool to ``reduce the number of variables'',
meaning starting  from a ``big'' Poisson manifold $(M,P)$, construct a
smaller one $(M_{red},P_{red})$.  Consider an
embedded coisotropic submanifold  in the Poisson manifold, $$\iota: C
\hookrightarrow M.$$
A submanifold of a Poisson manifold is called coisotropic iff
the vanishing ideal
$$
\mathcal{J}_C = \{f \in C^\infty(M) \; | \; \iota^*f = 0\} = \ker \iota^*.
$$
is closed under Poisson bracket.
This is equivalent to say that  $P^\sharp(N^*C)\subset TC$ \footnote{ On a Poisson manifold $(M,\PO)$, the map $\PO^\sharp$ is defined by $$\PO^\sharp: T^*M\rightarrow TM : \alpha \mapsto  \PO^\sharp(\alpha) \, \textrm{ so that }\, \beta(P^\sharp(\alpha):=P(\alpha,\beta).$$ }
where $N^*C(x)=\{\,\alpha_x\in T^*_xM\,\vert\, \alpha_x(X)=0\,\forall X\in T_xC\,\}$. 
In the symplectic case $P^\sharp(N^*C)=TC^{\perp}$ is the orthogonal with respect to the symplectic form $\omega$
of the tangent space to $C$ so that a submanifold $C$ in a symplectic manifold $(M,\omega)$ is isotropic iff $TC^{\perp}\subset TC$ .\\
 The characteristic distribution 
$P^\sharp(N^*C)$ is involutive; it is spanned at each point
by the Hamiltonian vector fields corresponding to functions\footnote{ The Hamiltonian vector field $X_f$ corresponding to a function $f\in \CiM$ is $X_f:= P^\sharp(df)$.} which are locally in $\mathcal{J}_C$.\\
We assume  the canonical foliation to have a nice leaf space $M_\red$, i.e. a structure of a
smooth manifold such that the canonical projection $\pi: C
\longrightarrow M_\red $ is a submersion.}
Then
$M_\red$ is a Poisson manifold in a canonical way:  defining
the normalizer $\mathcal{B}_C$ of $\mathcal{J}_C $
$$\mathcal{B}_C = \left\{f \in C^\infty(M) \; | \; \{f, \mathcal{J}_C\}
    \subseteq \mathcal{J}_C \right\},$$ 
    one has an isomorphism of spaces:
$$
  \mathcal{B}_C \big/ \mathcal{J}_C \simeq   \pi^* \Cinfty(M_\red) :
  [f] \; \mapsto \; \iota^*f 
  . 
$$
One defines the Poisson structure on $M_{red}$ to make it an isomorphism of Poisson algebras. 

Various procedure of reduction were presented in the context of deformation quantization.
Following  the reduction proposed by Bordemann, Herbig, Waldmann \cite{bib:BHW},
one starts from
the associative
algebra ${\deform{\mathcal{A}}} = (\Cinfty(M)[[\nu]], \star)$ which is
playing the role of the quantized observables of the big system. 
 A good analog of the vanishing ideal $\mathcal{J}_C$ will be a left
ideal ${\deform{\mathcal{J}}}_C \subseteq C^\infty(M)[[\nu]] $ such
that the quotient $C^\infty(M)[[\nu]] \big/ \deform{\mathcal{J}}_C$
is in $\mathbb{C}[[\nu]]$-linear bijection to the functions
$C^\infty(C)[[\nu]]$ on $C$. 
 Then one defines
$$ {\deform{\mathcal{B}}}_C = \{a \in {\deform{\mathcal{A}} } \; | \; [a,
{\deform{\mathcal{J}}}_C]_\star \subseteq {\deform{\mathcal{J}}}_C \},$$
 and one considers the associative
algebra $${\deform{\mathcal{B}}}_C \big/ {\deform{\mathcal{J}}}_C$$ as the
reduced algebra $\deform{\mathcal{A}}_{red}$. Clearly, 
one needs then to show that $\deform{\mathcal{B}}_C \big/
\deform{\mathcal{J}}_C$ is in $\mathbb{C}[[\nu]]$-linear bijection
to $\Cinfty(M_\red)[[\nu]]$ in such a way, that the isomorphism
induces a star product $\star_{red}$ on $M_\red$.  

We shall start from a
strongly invariant star product on $M$, and consider here the particular case of the Marsden-Weinstein
reduction: let $\mathsf{L}: G \times M \longrightarrow M$ be a smooth
left action of a connected Lie group $G$ on $M$ by Poisson
diffeomorphisms and assume we have an $\ad^*$-equivariant momentum
map $J$. The constraint manifold $C$ is  chosen to be the level surface
of $J$ for momentum $0 \in \g^*$ ( we assume, for simplicity,
that $0$ is a regular value) :  $ C = J^{-1}(\{0\})$; it is an embedded
submanifold which is coisotropic; we still denote by $\iota: C
\hookrightarrow M$ the inclusion. The group $G$ acts on $C$ and the
reduced space is the orbit space of this group action of $G$ on $C$.
In order to guarantee a good quotient we assume that $G$ acts freely
and properly and we  assume that $G$ acts properly not only on $C$ but on all of $M$. In
this case we can find an open $G$-invariant neighbourhood $M_{\mathrm{nice}} \subseteq M$ and  a $G$-equivariant
diffeomorphism
$$
    \Phi: M_{\mathrm{nice}} \longrightarrow U_{\mathrm{nice}} \subseteq C \times \g^*
$$
onto an open $g$-invariant neighbourhood $U_{\mathrm{nice}}$ of $C \times \{0\}$, where the
$G$-action on $C \times \g^*$ is the product action of the one on
$C$ and $\Ad^*$, such that for each $p \in C$ the subset $U_{\mathrm{nice}} \cap
(\{p\} \times \g^*)$ is star-shaped around the origin $\{p\}
\times \{0\}$ and the momentum map $J$ is given by the projection onto
the second factor, i.e. $J_{\vert_{M_{\mathrm{nice}}}}= \mathrm{pr}_2 \circ \Phi$.

BRST is a technique to describe the functions on the reduced space and was used in the theory of reduction in deformation quantization \cite{bib:BHW}; a simpler
description that we used with Waldmann \cite{bib:GW} is the classical Koszul resolution of $C^\infty(C)$.\\
The {\bf{Koszul complex}} is 
$
\left(\Cinfty(M_\nice, \AntiC^\bullet{\g}) = \Cinfty(M_\nice)
\otimes \AntiC^\bullet {\g}\,, \, \koszul\right)$,
  $\koszul$ being 
  the \emph{Koszul differential}  $
    \koszul x = \ins(J) x=\sum_a J_a \ins(e^a) x$
 with $    \{ e_a\} $  a basis of $\g.$\\
Defining the {\emph{prolongation map}} : $  \Cinfty(C)  \ni \phi \mapsto \prol(\phi) = (\pr_1 \circ \Phi)^* \phi \in
    \Cinfty(M_\nice)$, and
the \emph{homotopy} : $
\Cinfty(M_\nice, \AntiC^k {\g})\ni x \mapsto
 (h_k x)(p) = e_a \wedge \int_0^1 t^k 
    \frac{\partial (x \circ \Phi^{-1})}{\partial \mu_a}
    (c, t\mu)
    \D t,$
one shows that the Koszul complex is acyclic; also $
        \prol \iota^* + \koszul_1 h_0 = \id_{\Cinfty(M_\nice)}, \, \iota^* \koszul_1 = 0, \,
        h_0 \prol = 0,
$ and $$
    \Cinfty(M_\nice) \big/ (\mathcal{J}_C \cap \Cinfty(M_\nice))
    =  \ker \koszul_0 \big/ \Image \koszul_1
    \cong
    \Cinfty(C)
\,{\textrm{with}}\,\koszul_0=\iota^*.$$
If  $\mathcal{B}_C$ is the normalizer of $\mathcal{J}_C$,  the map :
\begin{equation*}
    \mathcal{B}_C \big/ \mathcal{J}_C \rightarrow   \pi^* \Cinfty(M_\red) :
  [f] \mapsto \iota^*f   
\end{equation*}
induces indeed an isomorphism of vector spaces  because 
 $\mathcal{J}_C = \Image \koszul_1$ and
$$f \in\mathcal{B}_C\,\textrm{ 
        iff}\, 0= \iota^*\{J_\xi,f\}=\iota^*(\Lie_{\xi^{*M}}f ) =\Lie_{\xi^{*C}}(\iota^*f) ~\forall\xi
        \in \lie{g}\,\textrm{ 
        iff}\, \iota^* f \in \pi^*
        \Cinfty(M_\red).$$
If $u \in \Cinfty(M_\red)$ then $\prol(\pi^*u) \in
\mathcal{B}_C$ so  the above map  is surjective;
the injectivity  is clear.\\
The Poisson bracket on $M_\red$ is  defined
through this bijection and gives explicitly
\begin{equation*}
    \pi^* \{u, v\}_\red = \iota^*\{\prol(\pi^* u), \prol(\pi^*v)\} \qquad u, v \in \Cinfty(M_\red).
\end{equation*}

Let $\star$ be a strongly invariant bidifferential formal star product\footnote{Recall that a star product is strongly invariant if it is invariant and
  $  J_\xi \star f - f \star J_\xi 
    = \nu \{J_\xi, f\}
    =  \nu \Lie_{\xi^{*M}} f
$
for all $f \in \Cinfty(M)[[\nu]]$ and $\xi \in \lie{g}$.} 
 on $M$, so that we start from the ``big" algebra of quantized observables
$$\deform{\mathcal{A}} = (C^\infty(M_\nice)[[\nu]], \star).$$
To define the left ideal, one first {\bf{deforms the Koszul complex}}, introducing \\ 
\emph{ a quantized Koszul operator}
$\deform{\koszul}^{(\kappa)}: \Cinfty(M_\nice, \AntiC^\bullet
    \lie{g})[[\nu]] \longrightarrow \Cinfty(M_\nice, \AntiC^{\bullet-1}
    \lie{g})[[\nu]]$ 
 defined by $$
        \deform{\koszul}^{(\kappa)} x
        =
        \ins(e^a) x \star J_a
        +
        \frac{\nu}{2} C_{ab}^c 
        e_c \wedge \ins(e^a) \ins(e^b) x
        +
       \nu \kappa \ins(\Delta) x,
$$
 where $C_{ab}^c = e^c([e_a, e_b])$ are the structure constants in the basis $\{e^a\}$, $\{e^a\}$ being the dual basis, and $\Delta\in \g^*$  is the modular form $ \Delta (\xi) = \tr \ad (\xi)$; 
 one checks that 
        $\deform{\koszul}^{(\kappa)}$ is left $\star$-linear,
         $\deform{\koszul}^{(\kappa)}$ is $G$-equivariant and 
          $\deform{\koszul}^{(\kappa)} \circ \deform{\koszul}^{(\kappa)}
        = 0$;\\
 \emph{a deformation of the restriction map :} $ \deformiotak   :\Cinfty(M_\nice)[[\nu]]
    \longrightarrow
    \Cinfty(C)[[\nu]]$ defined by 
    $$ \deformiotak=
    \iota^*\left(
        \id + \left(
            \deform{\koszul}^{(\kappa)}_1 - \koszul_1
        \right)
        h_0
    \right)^{-1}   
;$$ one can prove that although $h_0$ is not local, there exists  a formal series $S_\kappa = \id + \sum_{r=1}^\infty
    \lambda^r S^{(\kappa)}_r$ of $G$-invariant differential operators on $M_\nice$
    such that $\deformiotak = \iota^* \circ S_\kappa$ and  $S_\kappa 1 = 1$;\\
\emph{ a deformation of the homotopies;} 
   $ \deform{h}^{(\kappa)}_0 :\Cinfty(M_\nice)[[\nu]]
    \longrightarrow
    \Cinfty(M, \lie{g})[[\nu]]$ given by 
    $$ \deform{h}^{(\kappa)}_0=
    h_0 \left(
        \id + \left(
            \deform{\koszul}^{(\kappa)}_1 - \koszul_1
        \right)
        h_0
    \right)^{-1} \, \textrm{and higher terms}\,   \deform{h}^{(\kappa)}_k 
    = h_k 
    \left(
        h_{k-1} \deform{\koszul}^{(\kappa)}_k
        + \deform{\koszul}^{(\kappa)}_{k+1}  h_k
    \right)^{-1}.$$
All those maps are $G$-invariant, $ \deform{h}^{(\kappa)}_0 \prol = 0$, $ \deformiotak \deform{\koszul}^{(\kappa)}_1 = 0$, $\deformiotak \prol =
    \id_{\Cinfty(C)[[\nu]]}$, and \\
   $ \prol \deformiotak 
    + \deform{\koszul}^{(\kappa)}_1 \deform{h}^{(\kappa)}_0 =
    \id_{\Cinfty(M_\nice)[[\nu]]}$ as well as $ \deform{h}^{(\kappa)}_{k-1} \deform{\koszul}^{(\kappa)}_k
        +
        \deform{\koszul}^{(\kappa)}_{k+1} \deform{h}^{(\kappa)}_k
        =
        \id_{\Cinfty(M, \AntiC^k \lie{g})[[\nu]]}$.\\[1mm]       
One defines  the  {\bf{deformed left star ideal}}: 
$$ \deform{\mathcal{J}}_C = \image
    \deform{\koszul}^{(\kappa)}_1 =\ker  \deformiotak .
$$
 The left module
        $\Cinfty(M_\nice)[[\nu]] \big/ \deform{\mathcal{J}}_C$ is
        isomorphic to $\Cinfty(C)[[\nu]]$  with module structure
   $\bulletk$ defined by
    $
        f \bulletk \phi
        = \deformiotak (f \star \prol(\phi))$ for $\phi \in
    \Cinfty(C)[[\nu]],f \in \Cinfty(M_\nice)[[\nu]
$
via the map
      $$ \Cinfty(M_\nice)[[\nu]]  \big/ \deform{\mathcal{J}}_C\rightarrow \Cinfty(C)[[\nu]] : [f] \mapsto \deformiotak f \,\textrm{ whose inverse is
}\,   \phi \mapsto [\prol(\phi)] .$$
This left module structure
        is $G$-invariant ($
            \mathsf{L}^*_g  (f \bulletk \phi)
            = (\mathsf{L}^*_g f) \bulletk (\mathsf{L}^*_g \phi) $
        for all $g \in G$, $f \in \Cinfty(M)[[\nu]]$, and $\phi
        \in \Cinfty(C)[[\nu]]$) and for all $\xi \in
        \lie{g}$ one has, using the fact that the star product is strongly invariant,  
       $  J_\xi \bulletk \phi 
            = 
            \nu \Lie_{\xi^{*C}} \phi 
            - \nu \kappa \Delta(\xi) \phi.
$
\\[1mm]            
One considers  its {\bf{normalizer}} $$
    \deform{\mathcal{B}}_C
    =
    \left\{
        f \in \Cinfty(M)[[\nu]]
        \; \big| \;
        [f, \deform{\mathcal{J}}_C]_\star 
        \subseteq \deform{\mathcal{J}}_C
    \right\}
.$$
  For a $g$ in  $\deform{\mathcal{J}}_C$,
 $g = g^a \star J_a +
    \nu\kappa C_{ba}^a g^b$ with $g^a \in \Cinfty(M)[[\nu]]$;
 for $f \in
    \Cinfty(M)[[\nu]]$, the $\star$-bracket is 
  $
    [f, g]_\star 
    = \deform{\koszul}^{(\kappa)}_1 h 
    -\nu g^a \star \Lie_{(e_a)^{*M}} f
    $
    with  $h=(f\star g_a)e^a \in \Cinfty(M, \lie{g})$. Thus $[f,
    g]_\star$ is in $\deform{\mathcal{J}}_C$ iff $g^a \star
   \Lie_{(e_a)^{*M}} f$ is in the image of
    $\deform{\koszul}^{(\kappa)}_1$ for all $g^a$. This shows that $f
    \in \deform{\mathcal{B}}_C$ iff \\
    $\Lie_{\xi^{*M}} f \in \image
    \deform{\koszul}^{(\kappa)}_1 = \ker \deformiotak$
    thus
        iff $\Lie_{\xi^{*C}} \deformiotak f = 0$ for all $\xi
        \in \lie{g}$, i.e. iff $\deformiotak f \in \pi^*
        \Cinfty(M_\red)[[\nu]]$.\\[1.5mm] 
The {\bf{quotient algebra}}
        $\deform{\mathcal{B}}_C \big/ \deform{\mathcal{J}}_C$ is
        isomorphic to $\Cinfty(M_\red)[[\nu]]$ via the map        $$  \deform{\mathcal{B}}_C \big/ \deform{\mathcal{J}}_C\rightarrow \pi^*\Cinfty(M_\red)[[\nu]] :
             [f] \mapsto \deformiotak f 
       \, \textrm{whose inverse  is }\,        
             u \mapsto [\prol(\pi^*u)] .
        $$
        The 
         {\bf{reduced star}}
        product $\starredk$ on $\Cinfty(M_\red)[[\nu]]$ is induced from
        $\deform{\mathcal{B}}_C \big/ \deform{\mathcal{J}}_C$ and
        explicitly given by
        \begin{equation*}
            \label{eq:ReducedStarProduct}
            \pi^*(u \starredk v)
            = 
            \deformiotak \left(
                \prol(\pi^* u) \star \prol(\pi^* v)
            \right).
        \end{equation*}
   One checks that it is given by a series of bidifferential operators.\\
   
   In quantum mechanics, the algebra of quantum observables 
has a $^*$-involution given in the usual picture, where observables
are represented by operators, by the
passage to the adjoint operator. In  deformation
quantization,  a $^*$-involution on
$\deform{\mathcal{A}}=(C^\infty(M)[[\nu]], \star)$ for $\nu=i\lambda$ (with $\lambda\in\R)$ may be obtained, asking  the
star product to be Hermitian, i.e such that $\cc{f \star g} = \cc{g} \star \cc{f}$ and the $^*$-involution is then  complex
conjugation.  We have studied in \cite{bib:GW} how to
get in a natural way a $^*$-involution for the reduced algebra,
assuming that $\star$ is a Hermitian star product on $M$.  
The main idea here is to use a representation of the reduced quantum
algebra and to translate the notion of the adjoint. Observe that
$\deform{\mathcal{B}} \big/ \deform{\mathcal{J}}$ can be identified
(with the opposite algebra structure) to the algebra of
$\deform{\mathcal{A}}$-linear endomorphisms of $\deform{\mathcal{A}}
\big/ \deform{\mathcal{J}}.$ We use an additional positive
linear functional  $\omega: \deform{\mathcal{A}} \longrightarrow \mathbb{C}[[\lambda]]$
such that the
Gel'fand ideal of $\omega$, $\deform{\mathcal{J}}_\omega = \left\{a
    \in \deform{\mathcal{A}} \; \big| \; \omega(a^*a) = 0 \right\}$,
    coincides with the left ideal $\deform{\mathcal{J}}$ used in reduction, and 
     such that all
left $\deform{\mathcal{A}}$-linear endomorphisms of  the space of the GNS representation
$\mathcal{H}_\omega =\deform{ \mathcal{A}}
\big/ \deform{ \mathcal{J}}_\omega$, with the pre Hilbert space
structure defined via $\SP{\psi_a, \psi_b} = \omega(a^*b)$, are adjointable. 
Then the algebra of $\deform{\mathcal{A}}$-linear endomorphisms of
$\mathcal{H}_\omega$ (with the opposite structure) is equal to
$\deform{\mathcal{B}} \big/ \deform{\mathcal{J}}_\omega$ so that
$\deform{\mathcal{B}} \big/ \deform{\mathcal{J}}$ becomes in a natural
way a $^*$-subalgebra of the set $\Bounded(\mathcal{H}_\omega)$ of
adjointable maps.  

A
formal series of smooth densities $\sum_{r=0}^\infty \lambda^r \mu_r$
$\in \Gamma^\infty(|\Anti^{\mathrm{top}}| T^*C)[[\nu]]$ on the
coisotropic submanifold $C$, such that $\cc{\mu} = \mu$ is real, $\mu_0
> 0$ and so that $\mu$ transforms under the $G$-action as
$\mathsf{L}^*_{g^{-1}} \mu = \frac{1}{\Delta(g)} \mu$ (where $\Delta$
is the modular function), yields a positive linear functional which
defines a $^*$-involution on the reduced space. 
In the
classical Marsden Weinstein reduction, complex conjugation is a
$^*$-involution of the reduced quantum algebra. Looking  whether the $^*$-involution corresponding to a series of
densities $\mu$ is the complex conjugation yields  a new notion of
quantized modular class.

We also studied  in \cite{bib:GW} representations of the reduced algebra with the $^*$-involution given
by complex conjugation, relating the categories of modules
of the big algebra and the reduced algebra.  The usual technique to relate categories of modules is to use
a bimodule and the tensor product to pass from modules of one algebra
to modules of the other.  The construction of
the reduced star products gives a bimodule structure on
$\Cinfty(C)[[\nu]].$   
The space of formal series $\Cinftycf(C)[[\nu]]$ where
$$\Cinftycf(C) = \big\{ \phi \in \Cinfty(C) \; \big| \;
\supp(\phi) \cap \pi^{-1}(K) \; \textrm{is compact for all compact}K \subseteq M_\red\big\}$$  is a left
$(\Cinfty(M)[[\nu]], \star)$- and a right
$(\Cinfty(M_\red)[[\nu]], \starred)$-module; on  this bimodule there is a
$\Cinfty(M_\red)[[\nu]]$-valued inner product. This bimodule
structure and inner product on $\Cinftycf(C)[[\nu]]$ give a
strong Morita equivalence bimodule between
$\Cinfty(M_\red)[[\nu]]$ and the finite rank operators on
$\Cinftycf(C)[[\nu]]$.

 \section{Convergence of some formal star products}\label{section:convergence}
 
 For physics, $\hbar$ is a constant of nature and $\nu=i\hbar$ is not a formal parameter. Formal deformation  is not enough; for instance,
there is no general reasonable notion of spectra for formal star product algebras.  Spectra can be recovered only for a few examples  with convergence as  in \cite{bib:Bayen2}.
In general, formal deformation quantization can not predict the values one would obtain by  measurements.
In non formal deformation quantization of a Poisson manifold, one would like to have a subalgebra $\mathcal{A}$ of complex valued smooth functions (or distributions) on the manifold, with some topology,  and  a family of continuous associative law $*_\lambda$  on $\mathcal{A}$, depending on a parameter $\hbar$ belonging to a set $I$ admitting $0$ in its closure,  so that the limit of  $*_\hbar$ when $\hbar\mapsto 0$ is the usual product,  and the limit of the $\frac{[ \cdot, \cdot ]_{*_\hbar}}{\hbar}$ is the Poisson bracket. One would also like  the topology to be such that one could define nice representations of $(\mathcal{A},*_\lambda)$ and spectra. 
It is well known that the framework of $C^*$-algebras provides  a nice background for a  notion of spectra (the spectrum of an element $a$ in a unital $C^*$-algebra is the set of $\lambda\in\C$ such that $a-\lambda 1$ is not invertible), but this framework might be too restrictive.
 Formal deformation quantization  is  not  a solution but could be thought as a first step, using the constructions of that theory  to build, in a second step,  a framework where spectra and expectation values could be defined.  For a presentation of the convergence problem in deformation quantization, we recommend Waldmann's paper \cite{bib:Wconv}.\\

The Moyal star product presents interesting features concerning convergence. Recall that the formal Moyal star product comes  from  the quantization 
of  polynomials  on $\R^{2n}$ with Weyl's ordering.
 Weyl quantization can be extended beyond polynomials; heuristically one would like to write \\
$$``\mathcal{Q}_{Weyl} (F)"=(\frac{1}{2\pi})^{2n} \int_{\R^n}\int_{\R^n}\hat{F}(u,v)e^{i(uQ+vP)}dudv$$ where $\hat{F}$ is the Fourier transform
$\hat{F}(u,v)=\int_{\R^n}\int_{\R^n}{F}(q,p)e^{-i(uq+vp)}dqdp$.
If one develops  formally this, using the fact that
on a nice test function $\phi$, $(e^{iuQ}\phi)(x)=e^{iu.x}\phi(x)$,  $(e^{ivP}\phi)(x)=\phi(x+\hbar v)$ and $e^{i(uQ+vP)}=e^{\frac{i}{2}\hbar u.v}e^{iuQ}\circ e^{ivP}$,
 one gets the formula \\
 $(``\mathcal{Q}_{Weyl} (F)" (\phi))(x)=\frac{1}{(2\pi)^{2n}}\int\int \int \int_{(\R^n)^{\otimes4}}   F(q,p) e^{-isq-itp} e^{ i\hbar st/2} e^{ isx}
\phi(x+\hbar t)ds dt dq dp$.
If $t = \frac{y-x}{\hbar}$, we get
$ \frac{1}{(2\pi)^{2n} \hbar^n} \int\int \int \int_{(\R^n)^{\otimes4}} F(q,p) e^{-isq-i(y-x)p/\hbar}e^{ is(x+y)/2}
 f(y)ds dt dq dp$, which is 
$ \frac{1}{(2\pi \hbar)^n} \int_{\R^n} \int_{\R^n} F(\frac{x+y}{2},p) e^{-i(y-x)p/\hbar}\ dy dp$. Setting ${p = 2\pi \hbar \xi}$, it gives
$$
(\mathcal{Q}_{Weyl} (F) (\phi))(x):= \int_{\R^n}\left(\int_{\R^n} F(\frac{x+y}{2},2\pi\hbar\xi)e^{-2\pi i (y-x)\xi}\phi(y)dy\right)\,d\xi;
$$
which one takes as a definition of $\mathcal{Q}_{Weyl} (F)$; 
it is well defined for a test function $\phi$ in the Schwartz space when $F$ satisfies weak regularity bounds (there exists a constant $C>0$ and constants $C_{i,j}>0$ 
$\forall i,j \ge 0$ such that for all  $x,p$, one has $\vert\nabla^i_x\nabla^i_pF(x,p)\vert \le C_{i,j}(1+\vert x\vert +\vert p\vert)^C$). 
The above formula coincides with the previous one when $F$ is a  polynomial. \\
$\bullet$ The map $\mathcal{Q}_{Weyl}$ gives an isometry  between the space $L^2(\R^{2n})$ and the space of Hilbert Schmidt operators on $L^2(\R^n)$, associating a self-adjoint operator to a real function.\\
 $\bullet$ If $F$ and $G$ are two Schwartz functions, then the composition of the corresponding operators  $\mathcal{Q}_{Weyl} (F)\circ \mathcal{Q}_{Weyl} (G)$ is equal to $\mathcal{Q}_{Weyl}(F*_\hbar G)$ where $F*_\hbar G$ is the function defined by 
\begin{equation}
(F*_\hbar G)(u):
= (\frac{1}{\pi\hbar})^{2n} \int_{\R^{2n}} \int_{\R^{2n}} e^{\frac{2i}{\hbar}\Omega(v,w)} F(u+v)G(u+w)dvdw \label{compWeylgen}
\end{equation}
\begin{equation}
 = (\frac{1}{\pi\hbar})^{2n}\int_{\R^{2n}} \int_{\R^{2n}} e^{\frac{2i}{\hbar}(\Omega(u,v)+\Omega(v,w)+\Omega(w,u))} F(v)G(w)dvdw.    \label{compWeylgenS}
\end{equation}
with $\Omega={\tiny{\left( \begin{matrix} 0&I\cr -I&0\end{matrix}\right)}}$.
The result is a Schwartz function; hence $*_\hbar$ gives an associative product on the space of Schwartz functions, called the {\bf{convergent Moyal star product}}.\\
$\bullet$ The (formal) Moyal star product  can be seen as an asymptoptic expansion in $\nu = i\hbar$ of this composition law.\\

Many examples of star products are related to integral formulas.  For instance, the Berezin or Toeplitz star product on  
 K\"ahler manifolds  are obtained as asymptotic expansions for  $\hbar\rightarrow 0$  of some convergent counterpart in usual
quantization (see for instance \cite{bib:CGR} and  \cite{bib:BMS}), given by an integral formula. 
For instance, if $(M,\omega,J)$ is  a K\"ahler manifold and $(L,\nabla,h)$  is a
regular quantization bundle over $M$,  the formula for the composition of Berezin's symbols as defined in equations \eqref{eq:Bsymbol} and \eqref{comp} is given by 
\begin{equation}\label{eq:compsymb}
(A*_kB)(x)=\int_M\hat A(x,y)\hat B(y,x)\psi^k(x,y)\epsilon^{(k)}{k^n\omega^n\over
n!}\qquad A, B \in \ehat[l], k \geq l
\end{equation}
where $
\psi(x,y) = {|\langle e_{q'},e_q \rangle |^2\over \Vert e_{q'}\Vert^2 \Vert
e_q\Vert^2}$ with $\pi(q)=x$ and $\pi(q')=y.$ \\
The asymptotic expansion in
$k^{-1}$ as $k\to\infty$
is well defined; it gives a series in  $\frac{1}{k}$  which is  a differential star product on the manifold.

The difficulty to get convergent deformations in this framework of an integral formula depending on a parameter $k$ (given an associative law $*_k$ on a space $\ehat[k]$) is to find an algebra, i.e. a subspace stable by all $*_k$. \\
An interesting example is the disk; Berezin's procedure can be extended to non compact K\"ahler manifolds \cite{bib:CGR2}. For a possibly unbounded operator  ${\bf{A}}$ to have a Berezin's symbol, the coherent states must be in the domain, ${\bf{A}}e_q\in \H,\, \forall q$; do be able to write a composition formula ${\bf{A}}\circ {\bf{B}}$ in terms of symbols as above, one needs the adjoint of ${\bf{A}}$ to be defined on coherent states (so the section $s(x)=<e_{q'},{\bf{A}}e_q>q$ should be holomorphic and square integrable for all $q'$) and one needs  all ${\bf{B}}e_q$ to be in the domain of ${\bf{A}}$.\\ 
Consider the open disk, $\left(\Disk,\omega = {-i\lambda dz\wedge d\zbar\over2\pi(1-\vert z\vert^2)^2} = d\left({i\lambda\zbar
dz\over2\pi(1-\vert z\vert^2)}\right)\right)$; then $\Disk=SU(1,1)/U(1)$ and the
action of $SU(1,1)$ is Hamiltonian.
If $(L, \nabla, h)$ is a homogenous quantization for the simply-connected
group $\widetilde{SU(1,1)}$ then $L$ can be trivialised on all of $\Disk$ by
a section $s_0$ with
$
\vert s_0 \vert^2 = (1-\vert z\vert^2)^\lambda.
$
The norm on holomorphic sections is 
$$
\Vert fs_0 \Vert^2 = \int_{\Disk} \vert f(z)\vert^2 (1-\vert
z\vert^2)^\lambda {\lambda d^2z\over \pi(1-\vert z\vert^2)^2}
$$
where $d^2z$ denotes the usual Lebesgue measure;
$\Vert s_0\Vert^2$ is finite  for $\lambda >1$ which we assume. The characteristic function is 
$
\epsilon = 1 - \lambda^{-1}
$\\
The class of symbols which we shall use are the symbols of differential
operators $D(p,q,k)$ on $L^k$ defined by
$$
D(p,q,k)(fs^k_0)(z) = \left\{z^p \left({\partial\over\partial z}\right)^q
f(z)\right\}s^k_0(z).
$$
We have
$
D(p,q,k)e^{(k)}_{s_0(w)}(z) = \epsilon^{(k)} P_q(k\lambda)z^p
\left({\wbar\over 1-\wbar z}\right)^q(1-\wbar z)^{-k\lambda} s_0^k(z)
$,
where $P_q$ is the polynomial of degree $q$ given by $P_q(x) :=
x(x+1)\ldots(x+q-1)$, and the symbol of $D(p,q,k)$ is
given by
$$
\widehat{D(p,q,k)} (z) = P_q(k\lambda)\,z^p
\left({\zbar\over 1-\vert z\vert^2}\right)^q.
$$
It follows that $z^p\left({\zbar\over1-\vert z\vert^2}\right)^q$ is the
symbol of the densely defined operator ${D(p,q,k)\over P_q(k\lambda)}$
on $\H_k$. We can clearly compose such operators since the result of
applying the first to a coherent state is a coherent state for a
different parameter and these are in the domain of the second.
 So the $*_k$- defined in \eqref{eq:compsymb} is well-defined on those functions and yields
$$\aligned
\Biggl\{z^p\Biggl({\zbar\over1-\vert z\vert^2}\Biggr)^q\Biggr\} &*_k
\left\{z^r\left({\zbar\over1-\vert z\vert^2}\right)^s\right\}
= \left(P_q(k\lambda)P_s(k\lambda)\right)^{-1} \widehat{D(p,q,k)\circ
D(r,s,k)}\\
&= \sum_{m=0}^{\min(q,r)} \binom{q}{m} {r!\over(r-m)!}
{P_{s+q-m}(k\lambda)\over P_q(k\lambda)P_s(k\lambda)}
z^{p+r-m}\left({\zbar\over 1-\vert z\vert^2}\right)^{s+q-m}.
\endaligned$$
We deduce that $\left\{z^p\left({\zbar\over1-\vert
z\vert^2}\right)^q\right\} *_k \left\{z^r\left({\zbar\over1-\vert
z\vert^2}\right)^s\right\}$ is a rational function of $k$; hence
 the asymptotic expansion is
convergent on symbols of polynomial differential operators. \\
We have on the disk a subspace of smooth functions $\left\{z^p\left({\zbar\over1-\vert
z\vert^2}\right)^q\right\}$, with a family of associative products $\{*_k\}$. \\
The star product on the dual $\g^*$ of a Lie algebra $\g$ obtained via the bijection between polynomials on $\g^*$ and the universal enveloping algebra, has also an integral formula counterpart; so has the star product on the cotangent bundle of a Lie group.
Whether in general the asymptotics can be used  to recover the convergent quantization is a topic of  research. \\

In the framework of $C^*$-algebras, Rieffel introduced the notion of strict deformation quantization (see \cite{bib:Rieffel3,bib:Rieffel2,bib:Rieffel1}):
 A {\bf{strict deformation quantization}}  of a dense $*$-subalgebra $\mathbb{A}'$ of a $C^*$-algebra, in the direction of a Poisson bracket $\{.,.\}$ defined on $\mathbb{A}'$, is an open interval $I\subset \R$ containing $0$, and the assignment, for each $\hbar \in I$, of an associative product $\times_\hbar$, an involution $*_\hbar$ and a $C^*$-norm $\Vert ~\Vert_\hbar$ (for $\times_\hbar$ and $*_\hbar$) on $\mathbb{A}'$, which coincide for $\hbar=0$ to the original product, involution and $C^*$-norm on $\mathbb{A}'$, such that the corresponding field of $C^*$-algebras, with continuity structure given by the elements of $\mathbb{A}'$ as constant fields, is a continuous field of $C^*$-algebras, and such that for all $a,b \in\mathbb{A}'$, $\Vert \frac{(a \times_\hbar b-ab)}{i\hbar}-\{ a,b\}\Vert_\hbar \to 0$ as $\hbar\to 0$. A problem is that very few examples are known.\\
  Group actions appear here in an essential way : Rieffel  introduced a  general way to construct such $C^*$-algebraic deformations
based on a strongly continuous isometrical action of $\R^d$ on a $C^*$-algebra $\mathbb{A}$
$$
\alpha : \R^d\times \mathbb{A}\rightarrow \mathbb{A} : (x,a)\mapsto \alpha_xa.
$$
The product formula for the smooth vectors $\mathbb{A}^\infty$ with respect to this action is defined, using an oscillatory integral, choosing a fixed element $\theta$ in the orthogonal Lie algebra $so(d)$, by
$$
a\times_\hbar b:=a*^\alpha_\theta b :=(\frac{1}{\pi\hbar})^{d}\int_{\R^d\times \R^d} \alpha_x(a)\alpha_y(b)\exp^{\frac{2i}{\hbar} x\cdot \theta y}dxdy
$$
and it gives a pre $C^*$ associative algebra structure on $\mathbb{A}^\infty$. This
generalizes the  Weyl quantization of $\R^{2n}$. Indeed formula \eqref{compWeylgen} can be rewritten as
$$
F *_\hbar G=  (\frac{1}{\pi\hbar})^{2n}\int_{\R^{2n}\times\R^{2n}} \tau_v(F)\tau_w(G)\e^{{\frac{2i}{\hbar}}\Omega(v,w)} dvdw 
$$
where $\tau$ denotes the action of $\R^{2n}$ on functions on $\R^{2n}$ by translation.\\
Bieliavsky  and Gayral have generalized the construction
to  actions of Lie groups that admit negatively curved left-invariant K\"ahler structure. An important observation due to Weinstein is the relevance in the phase appearing in the product kernel (see equation \eqref{compWeylgenS}) of  the symplectic flux $S(x,y,z)={\small{\Omega(x,y)+\Omega(y,z)+\Omega(z,x)}}$ through a geodesic triangle that admits the points $x,y$ and $z$ as mid-points of its geodesic edges. This lead to the study of  symplectic groups which have a structure of symmetric symplectic spaces.  Bieliavsky  and his collaborators have  built,with increasing generality,  analogues of Weyl's quantization :  they
gave  universal deformation formulas for those groups and obtained  new examples of strict deformation quantization  \cite{bib:Biel1, bib:Biel2, bib:BiGa, bib:BiDS} .\\

A difficulty arising considering convergent star products given by integral formulas (like the convergent star product defined on the space of Schwartz functions on $\R^{2n}$ given by formula \eqref{compWeylgen}) is  to extend the construction to infinite dimensional cases, and such an extension is necessary  to have a  deformation quantization approach for quantum field theory.  \\

Another approach to the convergence problem is the following. Taking the formal power series defining a formal star product, one can ask for  convergence in a mathematically meaningful way. This  has been achieved  by Waldmann et al. in a growing number of  examples, for instance  the Wick star product on $\C^n$ and even in infinite dimension \cite{bib:wald2,bib:BRW}, the star product obtained by reduction on  the disk \cite{bib:beisWald,bib:KRSW} , the so-called Gutt star product on the dual of a Lie algebra \cite{bib:ESW2}, a Wick type star product on the sphere \cite{bib:WES}. 
They take a class of functions on which the star product obviously converges, build seminorms which garantee the continuity of the deformed multiplication, and extend the product  by continuity to the completion of the class . In this way,  they construct  topological non-commutative algebras, over $\C$ and not just over $C[[\nu]]$, essentially of Fr\'echet type.  They study Hilbert space representations of these algebras by   a priori unbounded operators
 \cite{bib:SW}.  A  nice short presentation of  results is given in \cite{bib:Wconv}.
Convergence of the Moyal star product on a Fr\'echet algebra had also been studied by Omori et al in \cite{bib:OMY4}.

\end{document}